\documentclass[11pt,notitlepage,twoside]{article}
\usepackage{latexsym}
\usepackage{amsmath}
\usepackage{amsfonts}
\usepackage{amssymb}
\renewcommand{\a }{\alpha }
\renewcommand{\b }{\beta }
\renewcommand{\d}{\delta }

\newcommand{\D }{\Delta }

\newcommand{\e }{\varepsilon }
\newcommand{\g }{\gamma}

\renewcommand{\l }{\lambda }

\newcommand{\n }{\nabla }
\newcommand{\var }{\varphi }

\newcommand{\Sig }{\Sigma}

\newcommand{\ov}{\overline}
\newcommand{\wtilde }{\widetilde}

\newcommand{\be}{\begin{equation}}
\newcommand{\ee}{\end{equation}}
\newenvironment{pf}{\noindent{\bf Proof.}\enspace}{%\rule{2mm}{2mm}
\hfill$\Box$\medskip}
\newenvironment{pfn}[1]{\noindent{\bf Proof of {#1}\enspace}}{%\rule{2mm}{2mm}
\hfill$\Box$\medskip}
\newcommand{\R}{\mathbb{R}}
\newcommand{\N}{\mathbb{N}}

\newtheorem{thm}{Theorem}[section]
\newtheorem{pro}[thm]{Proposition}
\newtheorem{lem}[thm]{Lemma}

\newtheorem{df}[thm]{Definition}
\numberwithin{equation}{section}

\author{{\large \textbf{Khalil EL MEHDI} }\\
{\it\small Facult\'e des Sciences et Techniques}\\
 {\it\small Universit\'e de Nouakchott, BP 5026}\\
 {\it\small Nouakchott, Mauritania}\\
{\it\small E-mail: khalil@univ-nkc.mr}\\
{\it\small and }\\
{\it\small The Abdus Salam ICTP, Trieste, Italy }
}

 \title{{\bf\Large{ On Conformal Paneitz 
  Curvature Equations in Higher Dimensional Spheres }}}
\begin{document}

\date{}
\maketitle
{\footnotesize
\noindent{\bf Abstract.}
We study the  problem of
prescribing the Paneitz curvature  on higher dimensional
 spheres. Particular attention is paid to the blow-up points, i.e. the critical points at infinity of the corresponding variational problem. Using  
topological tools and a careful analysis of the gradient flow lines in the neighborhood of such critical points at infinity, we prove some existence results.
}

\section{ Introduction }
Let $(M,g)$ be a smooth Riemannian manifold of dimension $n\geq 4$, with scalar curvature $R_g$ and Ricci curvature $Ric_g$.
In $1983$, Paneitz \cite{P} introduced in dimension four the following fourth order operator
$$
P_g^4  = \D _g^2  - \mbox{ div}_g(\frac{2}{3}R_g - 2 \mbox{Ric}_g ) \circ  d,
$$
where $div_g$ denotes the divergence and $d$ the de Rham differential.

This operator enjoys the analogous covariance property as the Laplacian in dimension two: under conformal change of metric  $\tilde{g} = e^{2u}g$ we have
$$
P_{\tilde{g}}^4 = e^{-4u}P_g^4.
$$

In \cite{Br1}, Branson generalized the Paneitz operator to $n$-dimensional Riemannian manifolds, $n\geq 5$. Such an operator is related to the Paneitz operator in dimension four in the same way the conformal Laplacian is related to the Laplacian in dimension two and is defined as:
$$
P^n_g = \D ^2_g - div_g(a_n S_g g + b_nRic_g ) \circ d + \frac{n-4}{2} Q^n_g,
$$
where
$$
a_n=\frac{(n-2)^2 +4}{2(n-1)(n-2)}, \qquad b_n =\frac{-4}{n-2}
$$
$$
Q^n_g = - \frac{1}{2(n-1)} \D _g S_g + \frac{n^3 -4n^2 + 16n-16}{8(n-1)^2(n-2)^2}S^2_g - \frac{2}{(n-2)^2} |Ric_g|^2.
$$
Under the conformal change of metric  $\tilde{g}=u^{4/(n-4)}g$, the conformal Paneitz operator enjoys the covariance property:
$$
P^n_g(u\var ) = u^{(n+4)/(n-4)} P^n_{\tilde{g}}(\var )\qquad \qquad \mbox{ for all } \var\in C^\infty (M),
$$
and the closely related fourth order curvature invariant $Q_g^n$ satisfies
\begin{eqnarray}\label{e:1}
P^n_g(u) =\frac{n-4}{2} Q^n_{\tilde{g}} u^{(n+4)/(n-4)}\quad\mbox{on}\quad M.
\end{eqnarray}
We call $Q_g^n$ the Paneitz curvature.
For more details about the properties of the Paneitz operator, see for example \cite{Br2}, \cite{BCY}, \cite{C}, \cite{CGY1}, \cite{CQY1}, \cite{CQY2}, \cite{CY1}, \cite{DHL}, \cite{G}, \cite{XY}.

A problem  naturally arises when looking at equation \eqref{e:1}: the problem of prescribing the Paneitz curvature, that is, given a smooth function  $ f: M \to \R$, does there exist a metric $\tilde{g}$ conformally equivalent to $g$ such that $Q^n_{\tilde{g}}=f$ ? From equation \eqref{e:1}, the problem is equivalent to finding a smooth solution $u$ of the equation
\begin{eqnarray}\label{e:2}
P^n_g(u) =\frac{n-4}{2} f u^{(n+4)/(n-4)}, \qquad u > 0\qquad \mbox{on}\qquad M.
\end{eqnarray}
The requirement about the positivity of $u$ is necessary for the metric $\tilde{g}$ to be Riemannian. Problem \eqref{e:2} is the analogue of the classical scalar curvature problem to which a wide range of activity has been devoted in the last decades (see for example the monograph \cite{A} and  references therein). On the other hand, to the author's knowledge, problem \eqref{e:2} has been studied in \cite{BE1}, \cite{BE2}, \cite{C}, \cite{DMO1}, \cite{DMO2} \cite{ER}, \cite{F}, \cite{XY} only.

In this paper, we are interested in the case where a noncompact group of conformal transformations acts on the equation so that Kazdan-Warner type conditions give rise to obstructions, as in the scalar curvature problem, see \cite{DHL} and \cite{WX}. The simplest situation is the following: let $(S^n,g)$ be the standard sphere, $n\geq 5$, endowed with its standard metric. In this case our problem is equivalent to finding a solution $u$ of the equation
\begin{eqnarray}\label{e:3}
 \mathcal{P}u:=\D^2u-c_n\D u+d_n u =K u^\frac{n+4}{n-4}, \qquad u > 0\qquad \mbox{ on } S^n,
\end{eqnarray}
where $c_n=\frac{1}{2}(n^2-2n-4)$, $d_n=\frac{n-4}{16}n(n^2-4)$ and where $K$ is a given function defined on $S^n$.

Our aim is to give sufficient conditions on $K$ such that problem \eqref{e:3} admits a solution. Our approach uses dynamical and topological methods involving the study of critical points at infinity of the associated variational problem, see Bahri \cite{B1}. Precisely, we extend the topological tools introduced by Bahri \cite{B2} to the framework of such higher order equations. Our method relies on the use of the  invariant introduced by Bahri \cite{B2}, which we extend to prove some existence results for problem \eqref{e:3}. To state our main results, we need to introduce the assumptions that we will use and some notations.\\
${\bf (A_1)}$\hskip 0.3cm We assume  that $K$ is a positive $C^3$-function on $S^n$ and which has only nondegenerate critical points $y_0,...,y_s$  with
$$
K(y_0)=\max K,\, -\D K(y_i)>0 \mbox{ for }i=0,1,\, -\D K(y_i)<0 \mbox{ for }i\geq 2
$$
and $ index(K,y_1)\ne n$.\\
Let $Z$ be a pseudo gradient of $K$ of Morse-Smale type, that is, 
the intersections of the unstable and stable manifolds of the
critical points of $K$ are transverse. We denote by $(n-k)$ the Morse index of $y_1$ and we set
\be \label{e:X}
X=\overline{W_s(y_1)},
\ee
where $W_s(y_1)$ is the stable manifold of $y_1$ for $Z$. Let us define
$$
B_2(X)=\{\a_1\d_{x_1}+\a_2\d_{x_2}/\a_i\geq 0, \a_1+\a_2=1,\,
x_i\in X\},
$$
 where $\d_x$ denotes the Dirac mass at $x$. For $a\in S^n$ and $\l>0$, let
$$
\wtilde{\d}_{(a,\l)}(x)=\frac{\b_n}{2^{\frac{n-4}{2}}}
\frac{\l^\frac{n-4}{2}}{\bigl(1+\frac{\l^2-1}{2}(1-\cos
 d(x,a))\bigr)^\frac{n-4}{2}},
$$
 where $d$ is the geodesic distance on $(S^n,g)$ and $\b_n= [(n-4)(n-2)n(n+2)]^{(n-4)/8}$.
After performing a stereographic projection $\Pi$ with the
 point $-a$ as pole, the function $\wtilde{\d}_{(a,\l)}$ is
 transformed into

 $$
\d_{(0,\l)}=\b_n\frac{\l^{\frac{n-4}{2}}}{(1+\l^2\mid
 y\mid^2)^{\frac{n-4}{2}}},
$$
 which is a solution of the problem (see \cite{Li})
$$
 \D^2u= u^\frac{n+4}{n-4} ,\, u>0\,\quad \mbox{ on } \quad \R^n.  
$$
We notice that problem \eqref{e:3} has a variational structure. The corresponding functional is
\be\label{functional}
J(u)=\left(\int_{S^n}
K|u|^{2n/(n-4)}\right)^{(4-n)/n}
\ee
defined on the unit sphere $\Sig$ of $H_2^2(S^n)$  equipped with the norm:
$$
\mid\mid u\mid\mid^2= \langle u,u\rangle_\mathcal{P}=\int_{S^n} \mathcal{P}u\cdot u = \int_{S^n} \mid \D u\mid^2 +c_n
\int_{S^n} \mid\n u\mid^2 + d_n\int_{S^n} u^2.
$$
We set $\Sig ^+ =\{u\in \Sig \mid u >0\}$ and for $\l$ large enough, we introduce a map
$f_\l : B_2(X)\to \Sig ^+$, defined by
$$
 (\a_1\d _{x_1}+\a_2\d_{x_2})\longrightarrow
\frac{\a_1 \tilde \d_{(x_1,\l )} + \a_2 \tilde\d _{(x_2,\l )}}
{||\a_1 \tilde\d_{(x_1,\l)}+\a_2 \tilde\d _{(x_2,\l )}||}.
$$
Then, $ B_2(X)$ and $f_\l ( B_2(X))$ are manifolds in dimension
$2k+1$, that is, their singularities arise in dimension $2k-1$ and
lower, see \cite{B2}. Recall that $k$ satisfies $k=n-index(K,y_1)$ and therefore the dimension of $X$ is equal to $k$.\\
Let $\nu^+$ be a tubular neighborhood of $X$ in $ S^n$.
We denote by $\nu^+(y)$, for $y\in X$, the fibre at $y$ of this
tubular neighborhood. For $\e_1>0$, $z_1,z_2\in X$ such that $z_1\neq z_2$ and  $-\D
K(z_i)>0$ for $i=1,2$, we introduce the following set
\begin{align*}
\Gamma_{\e_1}=\bigg\{ &
\sum_{i=1}^2 \frac{\tilde\d_{(z_i+h_i,\l_i)}}
{K(z_i+h_i)^{\frac{n-4}{8}}}+v \mid v\in H^2_2(S^n) \mbox{ satisfies
}(V_0), \\
 & || v-\ov{v}||<\e_1,\, \, \l_i>\e_1^{-1}\mbox{ for }i=1,2,\, h_i\in
\nu^+(z_i),\, \mid h_1\mid^2+\mid h_2\mid^2<\e_1\bigg\}\notag,
\end{align*}
where $\ov{v}$ is defined in Lemma \ref{l:22} (see below) and
where $(V_0)$ is the following conditions:
\begin{align}\label{V0}
(V_0):\quad
\langle& v,\varphi_i\rangle_\mathcal{P}=0\,\,\mbox{ for }\,\,  i=1,2\,\,\mbox{ and every }\\
&\varphi_i = \wtilde{\d}_{(a_i,\l_i)},\,\, \partial\wtilde{\d}_{(a_i,\l_i)}/\partial\l_i,\,\, \partial\wtilde{\d}_{(a_i,\l_i)}/(\partial a_i)_j,\, j=1,..., n,\notag\\
&\mbox{for some system of coordinates }\, (a_i)_1,..., (a_i)_{n}\,\,\mbox{on } S^n\,\,\mbox{near } a_i:= z_i+h_i.\notag    
\end{align}
We also assume that\\
${\bf (A_2)}$\hskip 0.3cm $z_1$ and $z_2$ are distinct of $y_0$, or if one is $y_0$, the other one is $y_1$.\\
 For $\d>0$ small,
the boundary of $\Gamma_{\e_1}$ (defined by $|| v-\ov{v}|| =\e_1$, or
$\l_1=\e_1^{-1}$,  or $\l_2=\e_1^{-1}$, or  $\mid h_1\mid^2+\mid h_2\mid^2=\e_1$) does not
intersect $J^{-1}(c_\infty(z_1,z_2)+\d)$, where
\be\label{016}
c_\infty(z_1,z_2)=\left(S_n \sum_{i=1}^2\frac{1}
{K(z_i)^{(n-4)/4}}\right)^{4/n}.
\ee
 We then set
\begin{eqnarray}
C_\d :=C_\d(z_1,z_2)=\Gamma_{\e_1}\cap J^{-1}(c_\infty(z_1,z_2)+\d).
\end{eqnarray}
For $\e_1$ and $\d$ small enough, $C_\d(z_1,z_2)$ is a closed Fredholm (noncompact) manifold without boundary of codimension $2k+2$.\\
 For $\l$
large enough, we define the intersection number (modulo 2) of
$W_u(f_\l(B_2(X)))$ with $C_\d(z_1,z_2)$ denoted by
\begin{eqnarray}
\tau(z_1,z_2)=W_u(f_\l(B_2(X))).C_\d(z_1,z_2) ,
\end{eqnarray}
where $W_u(f_\l(B_2(X)))$ is the unstable manifold of
$f_\l(B_2(X))$ for a decreasing pseudogradient $V$ for $J$ which
is transverse to $f_\l(B_2(X))$. Notice that the dimension of $W_u(f_\l(B_2(X)))$ is equal to $2k+2$ and the codimension of $C_\d (z_1,z_2)$ is equal to $2k+2$. Therefore, the number $\tau (z_1,z_2)$ is well defined (see \cite{M}). Our  main result is the following.
\begin{thm}\label{t:1}
Let $n\geq 9$. If $\tau(z_1,z_2)=1$ for a couple $(z_1,z_2)\in X^2$ satisfying $(A_2)$ and $-\D K(z_i)>0$ for $i=1,2$, then \eqref{e:3} has a solution.
\end{thm}

The aim of the next result is to give some conditions on the function $K$  which allow
us to have $\tau(z_1,z_2)=1$ for some couple $(z_1,z_2)$ and thus, we obtain a solution for \eqref{e:3} by Theorem \ref{t:1}.
 Let $z_1,z_2\in X$ be such that $-\D
K(z_i) >0$. We choose $\nu^+(z_i)$ such that
$K(z_i)=\max_{\nu^+(z_i)}K$ and $z_i$ is the unique critical point
of $K$ on $\nu^+(z_i)$.
\begin{thm}\label{t:2}
Let $n\geq 9$. There exist positive constants $C_0$, $C_1$
such that, if, for two points $z_1$ and $z_2$ of
$X$, the following conditions hold:
\begin{enumerate}
\item $w(z_1,z_2):=\frac{K(z_1)+K(z_2)}{2K(y_1)}-1\leq C_0$.
\item For some positive constant $\rho_0$,
\begin{align*}
& w^{\frac{n-6}{n-4}}(a_1,a_2)\bigl(\frac{1}{d(a_1,a_2)^2}+
\frac{1}{\rho_0^2}\bigr)+\frac{\mid\n K(a_i)\mid^2}{K(a_i)^2} +w^{1/2}(a_1,a_2)\frac{|D^2 K(a_i)|}{K(a_i)}\\
&\quad  + w^{1/3}(a_1,a_2)\sup_{B_{(a_i,\rho_0)}}\left(\frac{|D^3 K(x)|}{ K(a_i)}\right)^{2/3} \leq \frac{C_1}{1+(\frac{\sup K}{K(y_1)})^{\frac{n-4}{8}}}\left(\frac{-\D K(a_i)}{K(a_i)}\right)
\end{align*}
 for each $i=1,2$, and for each $(a_1,a_2)\in \nu^+(z_1)\times
\nu^+(z_2)$ such that\\ $c_\infty(a_1,a_2)\leq c_\infty(y_1,y_1)$.
 \item $\displaystyle{\inf_{\partial(\nu^+(z_1)\times\nu^+(z_2))}
c_\infty(a_1,a_2)\geq c_\infty(y_1,y_1)}$,
\end{enumerate}
then \eqref{e:3} has a solution.
(Here $c_\infty(a_1,a_2)$ (resp $c_\infty(y_1,y_1)$ is defined by \eqref{016} replacing $(z_1,z_2)$ by $(a_1,a_2)$ (resp $(y_1,y_1)$)).
\end{thm}

The rest of the present paper is organized as follows. In
Section 2, we recall some preliminaries, introduce some definitions and the notations needed in the proof of our results.
In Section 3, we set up the variational structure and we perform an expansion of the Euler functional
associated to \eqref{e:3} and its gradient near the potential critical
points  at infinity. Then, we characterize the critical points at
infinity in Section 4. Lastly, Section 5   is
devoted to the proof of our results.

\section{ Preliminaries }
Solutions of problem \eqref{e:3} correspond, up to some positive constant, to critical points of the following functional defined on the unit sphere of $H^2_2(S^n)$ by
$$
J(u)=  \left(\int_{S^n}K|u|^{\frac{2n}{n-4}}\right)^{\frac
{4-n}{n}}.
$$
The exponent $2n/(n-4)$ is critical for the Sobolev embedding $H^2_2(S^n) \hookrightarrow L^q(S^n)$. As this embedding is not compact, the functional $J$ does not satisfy the Palais-Smale condition and therefore standard variational methods cannot be applied to find critical points of $J$. In order to describe the sequences
failing the Palais-Smale condition, we need to introduce some notations.
For $p\in \N^*$ and $\e>0$, we set
\begin{align*}
V(p&,\e) = \biggl\{u\in  \Sig \mid \exists a_1, ..., a_p \in S^n, \exists \l _1,..., \l _p > \e^{-1}, \exists \a _1,..., \a _p > 0  \mbox{ with }\\ 
 &\bigg |\bigg | u-\sum_{i=1}^p\a_i\wtilde{\d}_{(a_i, \l_i)}
\bigg |\bigg | < \e, \, \,  \e _{ij} < \e \, \forall i\neq j,\,  
\bigg |J(u)^{\frac{n}{n-4}}\a _i ^{\frac{8}{n-4}}K(a_i)-1\bigg |<\e \, \forall i  \biggr\},
\end{align*}
where
$$
\e_{ij}=\left(\frac{\l_i}{\l_j}+\frac{\l_j}{\l_i}+\frac{\l_i\l_j}{2}(1-\cos d(a_i,a_j))\right)^{(4-n)/2}.
$$
Let $w$ be a nondegenerate solution of \eqref{e:3}. We also set
$$
V(p,\e,w)=\biggl\{u\in \Sig |\, \exists \, \a_0>0\mbox{ with }
(u-\a_0w)\in V(p,\e) \mbox{ and } |\a_0 J(u)^{n/8}-1|<\e\biggr\}
$$

The failure of the Palais-Smale condition can be described,
following the ideas introduced in \cite{BrC}, \cite{L}, \cite{S},
as follows:
\begin{pro}\label{p:21}
Let $(u_j)\in \Sig^+$ be a sequence such that $\n J(u_j)$ tends to
zero and $J(u_j)$ is bounded. Then, there exist an integer $p\in
\N^*$, a sequence $\e_j>0$, $\e_j$ tends to zero, and an extracted
sequence of $u_j$'s, again denoted $u_j$, such that $u_j\in
V(p,\e_j,w)$ where $w$ is zero or a solution of \eqref{e:3}.
\end{pro}

The following lemma defines a parametrization of the set $V(p,\e)$. It follows from the corresponding statements in \cite{B2} and \cite{BC}.
\begin{lem}\label{l:21}
For any $p\in\N^*$, there is $\e_p>0$ such that if $\e\leq \e_p$ and $u\in
V(p,\e)$, then the following minimization problem
$$
\min\bigg\{\bigg|\bigg|u-\sum_{i=1}^p\a_i\tilde\d_{(a_i,\l_i)}
\bigg|\bigg|,\, \a_i>0,\, \l_i>0,\, a_i\in  S^n \bigg\}
$$ has a unique solution
$(\a,\l,a)=(\a_1,...,\a_p,\l_1,...,\l_p,a_1,...,a_p)$. In particular, we can write $u$ as
follows:
$$
u=\sum_{i=1}^p\a_i\tilde\d_{(a_i,\l_i)}+v,
$$
where $v$ belongs to $H^2_2(S^n)$ and satisfies
$(V_0)$.
\end{lem}

Next, we recall the following result which deals with the $v$-part of $u$.
\begin{lem}\label{l:22}\cite{BE1}
Assuming  the $\e_{ij}$'s are small enough and $J(u)^{\frac{n}{n-4}}\a_r^{\frac{8}{n-4}}K(a_r)$ is close to $1$  for $i\ne j$ and for $r=i,j$ , then there exists a unique $\ov{v}=\ov{v}(a,\a,\l)$ which minimizes \\$J\left(\sum_{i=1}^p \a_i \tilde{\d}_{(a_i,\l_i)}+v\right)$
with respect to $v\in E_\e:=\{v\mid v \,\mbox{satisfies}\,\, (V_0)\,\,\mbox{and}$ $ \mid\mid v\mid\mid < \e\}$, where $\e$ is a fixed small positive constant depending only on $p$. Moreover, we have the following estimate
\begin{eqnarray*}
\mid\mid \ov{v}\mid\mid\leq c \biggl[ 
\sum_{i=1}^p\left(\frac{\mid\n K(a_i)\mid}{\l_i}+\frac{1}{\l_i
  ^2}\right) + \sum_{i\ne
  j}\e_{ij}^{\min \left(1, \frac{n+4}{2(n-4)}\right)}(\log\e_{ij}^{-1})^{\min\left(\frac{n-4}{n}, \frac{n+4}{2n}\right)}\biggr].
\end{eqnarray*}
\end{lem}

Note that Lemma \ref{l:21} extends to the more general situation where the sequence $(u_j)$ of $\Sigma^+$, described in Proposition \ref{p:21}, has a nonzero weak limit, a situation which might occur if $K$ is the Paneitz curvature (up to a positive constant) of a metric conformal to  the standard metric $g$. Notice  that such a weak limit is a solution of \eqref{e:3}. Denoting by $w$  a nondegenerate solution of \eqref{e:3}, we then have the following lemma which follows  from the corresponding statement in \cite{B2}.
\begin{lem}\label{l:23}
For any $p\in\N^*$, there is $\e_p>0$ such that if $\e\leq \e_p$ and $u\in
V(p,\e,w)$, then the following minimization problem
$$\min\bigg\{\bigg|\bigg|u-\sum_{i=1}^p\a_i\tilde\d_{(a_i,\l_i)}-\a_0(w+h)
\bigg|\bigg|,\, \a_i>0,\, \l_i>0,\, a_i\in  S^n, \, h\in
T_w(W_u(w)) \bigg\}$$ has a unique solution
$(\a,\l,a, h )=(\a_1,...,\a_p,\l_1,...,\l_p,a_1,...,a_p,h)$. In particular, we can write $u$ as follows:
$$u=\sum_{i=1}^p\a_i\tilde\d_{(a_i,\l_i)}+\a_0(w+h)+v,$$
where $v$ belongs to $H^2_2(S^n)\cap T_w(W_s(w))$ and satisfies
$(W_0)$. Here $T_w(W_u(w))$ and $T_w(W_s(w))$ denote the tangent
spaces at $w$ of the unstable and stable manifolds of $w$, and $(W_0)$ are the following conditions:
\begin{equation}\label{W0}
(W_0):\quad
\begin{cases}
\langle v,\varphi_i\rangle_\mathcal{P}=0\,\,\mbox{ for }\,\,  i=1,...,p\,\,\mbox{ and every }\notag\\
\varphi_i = \wtilde{\d}_{(a_i,\l_i)},\,\, \partial\wtilde{\d}_{(a_i,\l_i)}/\partial\l_i,\,\, \partial\wtilde{\d}_{(a_i,\l_i)}/\partial (a_i)_j,\, j=1,..., n,\notag\\
\mbox{for some system of coordinates }\, (a_i)_1,..., (a_i)_{n}\,\,\mbox{on } S^n\,\,\mbox{near } a_i,\notag \\
\langle v, w\rangle =0,\notag\\
 \langle v, h_1\rangle =0  \quad\forall h_1\in T_w(W_u(w)).
\end{cases}
\end{equation}
\end{lem}

Now, following Bahri \cite{B2}, we introduce the following definitions and notations. 
\begin{df}\label{d:21}
A critical point at infinity of $J$ on $\Sigma^+$ is a limit of a flow-line $u(s)$  of equation $\frac{\partial u}{\partial s}= - \n J(u)$ with initial data $u_0\in\Sigma^+$ such that $u(s)$ remains in $V(p,\e(s),w)$ for large $s$. Here $w$ is zero or a solution of \eqref{e:3}, $p\in\N^*$, and $\e(s)$ is some function such that $\e(s)$ tends to zero when the flow parameter $s$ tends to $+\infty$. By Lemma \ref{l:23}, we can write such $u(s)$ as
$$
u(s)=\sum_{i=1}^p\a_i(s)\tilde\d_{(a_i(s),\l_i(s))}+\a_0(s)(w+h(s))+v(s).
$$
 Denoting $a_i=\lim_{s\to +\infty} a_i(s)$, we call $(a_1,...,a_p,w)_\infty$ a critical point at infinity of $J$. If $w\ne 0$, $(a_1,...,a_p,w)_\infty$ is called a mixed type of critical points at infinity of $J$.
\end{df}

In the sequel, we denote by  $A$  the set of $w$ such that $w$ is a critical point or a critical point at infinity of $J$ in $\Sig^+$ not containing $y_0$ in its description. We also denote by  $A_q$  the subset of $A$ such that the Morse
index of the critical point (at infinity) is equal to $q$. 
\begin{df} ({\it A family of pseudogradients $\mathcal{F}$}) A decreasing pseudogradient $V$ for $J$ is said to belong to
$\mathcal{F}$ if the following properties hold:\\
 - the set of  critical points at infinity of $J$ on $\Sig^+$ does not change if we take $V$ instead of $-\n J$ in the definition \ref{d:21},\\
- $V$ is transverse to $f_\l(B_2(X))$,\\
- for any $w\in A$, $(y_0,w)_\infty$ is a
critical point at
infinity with the following property:
\begin{align*}
 i((y_0,w)_\infty,w) &=1  & \forall &w \in A\\
 i((y_0,w)_\infty,w') &=0 & \forall &w'\in A,\, w'\ne w,\,
 \mbox{ index}(w')=\mbox{ index}(w)\\
 i((y_0,w)_\infty,(y_0,w')_\infty) &=i(w,w') & \forall &w'\in A,\,
 \mbox{ index}(w')=\mbox{ index}(w)-1.
 \end{align*}
Here and below $i(\var_1,\var_2)$ denotes the intersection number for $V$ of  $\var_1$ and $\var_2$ (see \cite{M}
and \cite{B2}) where $\var_i$ is any  critical point or a 
critical point at infinity of $J$.
 \end{df}
 \begin{df}\label{d:o}
Given a decreasing pseudogradient $V$ for $J$. We denote by
$\varphi(s,.)$ the associated flow. A critical point at
infinity $z_\infty$ is said to be dominated by $f_\l(B_2(X))$ if
$$\ov{\cup_{s\geq 0}\varphi(s,f_\l(B_2(X)))}\cap W_s(z_\infty)\ne
\emptyset.$$
\end{df}
Near the critical points at infinity, a Morse Lemma can be
completed (see Proposition \ref{c:310} and \eqref{mokh} below) so that the usual Morse theory can be extended and the
intersection can be assumed to be transverse. Thus the above
condition is equivalent to (see Proposition 7.24 and Theorem 8.2
of \cite{BR})
$$\cup_{s\geq 0}\varphi(s,f_\l(B_2(X)))\cap W_s(z_\infty)\ne
\emptyset.$$
\begin{df}
$z_\infty$ is said to be dominated by another critical point at
infinity $z'_\infty$ if
$$W_u(z'_\infty)\cap W_s(z_\infty)\ne \emptyset.$$
\end{df}

If we assume that the intersection is transverse, then
$index(z'_\infty)\geq index(z_\infty)+1$.

Given $w_{2k+1}\in A_{2k+1}$ and $V\in \mathcal{F}$, we denote by
\be\label{d:3}
(y_0,w_{2k+1})_\infty.C_{\d}
\ee
 the intersection number (modulo
2) of $W_u((y_0,w_{2k+1})_\infty)$ and $C_\d$.

In order to compute this intersection number, one can perturb $V$
(not necessarily in $\mathcal{F}$) so as to bring
$W_u((y_0,w_{2k+1})_\infty)\cap C_\d$ to be transverse. This
number is the same for all such small perturbations (just as in
degree theory). Notice that the dimension of $W_u((y_0,w_{2k+1})_\infty)$ is equal to $2k+2$ and the codimension of $C_\d$ is $2k+2$. Then $(y_0,w_{2k+1})_\infty.C_{\d}$ is also well defined, because the closure of $W_u((y_0,w_{2k+1})_\infty)$ only adds to $W_u((y_0,w_{2k+1})_\infty)$ the unstable manifolds of
critical points of index less than or equal to $2k+1$. These
manifolds are then of dimension $2k+1$ at most. Since the
codimension of $C_\d$ is equal to $2k+2$, these manifolds can be
assumed to avoid $C_\d$.

Now, for  $w_{2k+1}\in A_{2k+1}$ and $V\in \mathcal{F}$, we denote by
\be\label{d:27}
f_\l(B_2(X)).w_{2k+1}:=f_\l(B_2(X)).W_s(w_{2k+1})
\ee
the intersection number of $f_\l(B_2(X))$ and $W_s(w_{2k+1})$. We notice that the dimension of \\$f_\l (B_2(X))$ is equal to $2k+1$ and the codimension of $W_s(w_{2k+1})$ is equal to $2k+1$. Then, the intersection number, defined in \eqref{d:27} is well defined because $V$ is transverse
to $f_\l(B_2(X))$ outside $f_\l(B_1(X))$, which cannot dominate
critical points of index $2k+1$. Furthermore, $\ov{W_s(w_{2k+1})}$
adds to $W_s(w_{2k+1})$ stable manifolds of critical points of
 an index larger than or equal to $2k+2$. Since $f_\l(B_2(X))$ is of
dimension $2k+1$, these manifolds can be assumed to avoid it.

Lastly, we set for each $V\in\mathcal{F}$
\be\label{d:28}
I(V)=\tau - \sum_{w_{2k+1}\in A_{2k+1}}((y_0,
w_{2k+1})_\infty.C_\d)(f_\l(B_2(X)).w_{2k+1}).
\ee
Notice that \ref{d:28} was introduced by Bahri in \cite{B2} where he proved that $I(V)$ is independent on $V\in\mathcal{F}$. Namely, he showed in \cite{B2} that $I(V)=0,$ for each $V\in\mathcal{F}$ for the scalar curvature problem on $S^n$ with $n\geq 7$. We will prove that the same holds for the Paneitz curvature equation when $n\geq 9$.

\section{Expansion of the functional and its gradient }
This section is devoted to a useful expansion of $J$ and its
gradient near a critical point at infinity. In
order to simplify the notations, in the remainder we write
$\tilde\d_i$ instead of $\tilde\d_{(a_i,\l_i)}$. First, we prove the following result:
\begin{pro}\label{p:35}
For $\e >0$ small and $u=\sum_{i=1}^p\a
_i\tilde\d _{(a_i,\l _i)}+\a_0(w+h)+v\in V(p,\e,w )$, the following expansion holds
\begin{align*}
J(u) = & \frac{S_n\sum_{i=1}^p\a_i^2+\a_0^2||w||^2}{
(S_n\sum_{i=1}^p\a _i
^{\frac{2n}{n-4}}K(a_i)+\a_0^{\frac{2n}{n-4}}||w||^2
)^{\frac{n-4}{n}}} \left[1-\frac{c_2(n-4)}{n\b_0}
\sum_{i=1}^p\a _i ^{\frac{2n}{n-4}}\frac{4\D K(a_i) }{\l _i^2}\right.\\
 & -\frac{c_1}{\gamma_0}\sum_{i\ne j\geq 1}\a _i\a _j\e _{ij}
 +\frac{1}{\gamma_0}\left(Q_1(v,v)-f_1(v)\right)+\frac{\a_0^2}{\gamma_0}\left(Q_2(h,h)+f_2(h)\right)
\end{align*}
$$
 \left.+o\left(\sum_{i\ne j \geq 1}\e_{ij}+\sum_{i=1}^p\frac{1}{\l _i^2}+||v||^2+||h||^2\right)\right]
$$
%\end{align*}
where $c_1=\b_n^{2n/(n-4)}\int_{\R^n} \frac{dx}{(1+|x|^2)^{(n+4)/2}}$, $c_2 = \frac{1}{2n}\int_{\R^n}|x|^2 \d_{(0,1)}^{2n/(n-4)}$,\\ $S_n=\int_{R^n}\d_{(0,1)}^{2n/(n-4)}$, and where
\begin{align*}
Q_1(v,v) &= ||v||^2-\frac{n+4}{n-4}\left(\sum_{i=1}^p  \int_{S^n}\tilde\d _i^{\frac{8}{n-4}}
v^2+\int_{S^n}Kw^{\frac{8}{n-4}}v^2\right), \\
Q_2(h,h)&= ||h||^2-\frac{n+4}{n-4}
\int_{S^n}Kw^{8/(n-4)}h^2,\\
 f_1(v) &= \frac{2\gamma_0}{\b_0 }\int_{S^n}K\biggl(\sum_{i=1}^p\a
_i\tilde\d _i\biggr)^{(n+4)/(n-4)}v,\\
f_2(h)&=\frac{1}{\a_0}\sum_i\a_i\langle \tilde\d_i,h\rangle _{\mathcal{P}}-\frac{2\g_0}{\a_0\b_0}
\int_{S^n}K\biggl(\sum_{i=1}^p\a_i\tilde\d_i\biggr)^{(n+4)/(n-4)}h,\\
 \b_0  &=  S_n(\sum_{i=1}^p\a
_i^{2n/(n-4)}K(a_i))+\a_0^{2n/(n-4)}||w||^2, \\
\gamma_0 &= S_n(\sum_{i=1}^p\a _i^2)+\a_0^2||w||^2.
\end{align*}
\end{pro}
\begin{pf}
We recall that we have $\langle v,w \rangle_\mathcal{P} = \langle v,h\rangle_\mathcal{P} = \langle v,\tilde \delta_i\rangle_\mathcal{P} = \langle w,h\rangle_\mathcal{P}=0$.
  We need to estimate
$$
 N(u)= ||u||^2 \mbox{ and } D=\int_{S^n} K(x)
u^{\frac{2n}{n-4}}.
$$
We  have
$$
N(u)=\sum_{i=1}^p\a_i^2||\tilde\d_i||^2+\a_0^2(||h||^2+||w||^2)
+||v||^2+\sum_{i\ne j}\a_i\a_j\langle \tilde\d _i,\tilde\d_j\rangle +
2\sum_{i=1}^p\a_i\a_0\langle \tilde\d _i,w+h\rangle.
$$
Observe that
\begin{align*}
||\tilde\d_i ||^2 & =\int_{\R^n}|\D \d_i |^2 = S_n,\\
\langle\tilde\d _i,\tilde\d _j\rangle_\mathcal{P} & =\int_{\R^n}\d _i^{(n+4)/(n-4)} \d
_j= c_1\e _{ij} + O(\e _{ij}^{n/(n-4)}\log(\e
_{ij}^{-1})),\\
\langle\tilde\d_i,w\rangle_\mathcal{P} & =\int_{S^n}\tilde\d_i
^{(n+4)/(n-4)}w=O\left(\l_i^{(4-n)/2}\right).
\end{align*}
Thus
\begin{eqnarray}\label{N}
N = & \gamma_0 + c_1\sum_{i\ne j}\a _i\a _j\e _{ij}
+\a_0^2||h||^2+ ||v||^2+\a_0\sum_{i}\a_i\langle\tilde\d_i,h\rangle_{\mathcal{P}}\end{eqnarray}
$$+
o\left(\sum_{i=1}^p\frac{1}{\l_i^2}+ \sum_{i\ne j}\e _{ij}\right).
$$
%\end{align}
For the denominator, we write
\begin{align*}
D&=\int K (\sum_{i=1}^p\a _i\tilde\d _i
)^{\frac{2n}{n-4}}+\a_0^{\frac{2n}{n-4}}\int K
(w+h)^{\frac{2n}{n-4}}\\
&+ \frac{2n}{n-4} \int K
(\sum_{i=1}^p\a _i\tilde
\d _i+\a_0(w+h))^{\frac{n+4}{n-4}}v
  +\frac{2n\a_0}{n-4}\int K
(\sum_{i=1}^p\a _i\tilde \d _i)^{\frac{n+4}{n-4}}(w+h)\\
& +\frac{n(n+4)}{(n-4)^2}\int K(\sum_{i=1}^p\a
_i\tilde \d _i+\a_0(w+h))^{\frac{8}{n-4}}v^2 + O\left(\sum_{i=1}^p\int \tilde
\d_i(w+h)^{\frac{n+4}{n-4}}\right)\\
 & +O\left(\int_{S^n}(\sum_{i=1}^p \a_i\tilde\d_i)^{\frac{8}{n-4}}\min{}^2
(\sum_{i=1}^p\a_i\tilde\d_i,w+h)\right)+O\left(||v||^{\min(3,\frac{2n}{n-4})}\right)
\end{align*}
Observe that
\begin{align}
\int_{S^n}K (\sum_{i=1}^p\a _i\tilde\d _i)^{\frac{2n}{n-4}} & =
\sum_{i=1}^p\a_i ^{\frac{2n}{n-4}}\biggl(K(a_i)S_n+
c_2\frac{4\D K(a_i)}{\l_i^2}\biggr)\\
 & +\frac{2n}{n-4}\sum_{i\ne j}\a_i ^{\frac{n+4}{n-4}}\a_j K(a_i)
c_1\e _{ij}+o\left(\sum \e _{ij}+\sum \frac{1}{\l_i^2}\right).\notag
\end{align}
Using the fact that $h$ belongs to the tangent space at $w$, we
derive that
\begin{align}
\int_{S^n}K (w +h)^{\frac{2n}{n-4}}&  =  \int_{S^n}K
w^{\frac{2n}{n-4}}+\frac{2n}{n-4}\int_{S^n}K
w^{\frac{n+4}{n-4}}h+\frac{n(n+4)}{(n-4)^2}\int_{S^n}K
w^{\frac{8}{n-4}}h^2\notag\\
 &\, \,  +O(||h||^{\min(3,\frac{2n}{n-4})})\notag\\
 & =||w||^2+\frac{n(n+4)}{(n-4)^2}\int_{S^n}K
w^{\frac{8}{n-4}}h^2+O(||h||^{\min(3,\frac{2n}{n-4})}).
\end{align}
Since $v\in T_w(W_s(w))$ and $h\in T_w(W_u(w))$, the linear form
on $v$ can be written as
\begin{align}
&\int_{S^n}K (\sum_{i=1}^p\a _i\tilde \d _i 
+\a_0(w+h))^{\frac{n+4}{n-4}}v =\int_{S^n}K (\sum_{i=1}^p\a
_i\tilde \d _i)^{\frac{n+4}{n-4}}v\notag\\
&\quad+\int_{S^n}K
(\a_0(w+h))^{\frac{n+4}{n-4}}v
  +O\biggl(\sum_{i=1}^p\int\tilde\d_i
^{\frac{8}{n-4}}|w+h||v|+\int\tilde\d_i|w+h|^{\frac{8}{n-4}}|v|\biggr)\notag\\
 &\quad =\frac{\b_0}{2\gamma_0}f_1(v)+\a_0^{\frac{n+4}{n-4}}\left(\int
 Kw^{\frac{n+4}{n-4}}v+{\frac{n+4}{n-4}}\int
 Kw^{\frac{8}{n-4}}hv\right)\notag\\
&\quad\quad +O\biggl(||v||||h||^{\min(2,{\frac{n+4}{n-4}})}\biggr)
 \notag \\
 &\quad  =\frac{\b_0}{2\g_0}f_1(v)+O\biggl(||v||^{\min(3,{\frac{2n}{n-4}})}+
 ||h||^{\min(3,{\frac{2n}{n-4}})}\biggr).
\end{align}
 Furthermore, we have
\begin{align}\label{v2}
\int K(\sum_{i=1}^p\a _i\tilde \d
_i+\a_0(w+h))^{\frac{8}{n-4}}v^2= & \sum_{i=1}^pK(a_i)\int (\a _i\tilde \d
_i)^{\frac{8}{n-4}}v^2 \\
 & +\int K(\a _0 w)^{\frac{8}{n-4}}v^2
+o(||v||^2+||h||^2).\notag
\end{align}
Finally, we notice that
\be\label{ee:}
\langle \tilde\d_i, h\rangle  \int K\left(\sum\a_i\tilde\d_i\right)^{\frac{n+4}{n-4}}h=o\left(||h||^2\right); \,\langle \tilde\d_i, h\rangle f_1(v)=o\left(||h||^2 + ||v||^2\right).
\ee
Combining \eqref{N},...,\eqref{ee:} and the fact that
\be\label{ed}
J(u)^{n/(n-4)}\a_i^{{8}/(n-4)}K(a_i)=1+o(1)\,\forall i; \quad \a_0 J(u)^{n/8}=1+o(1),
\ee
the result follows.\end{pf}
\begin{pro}\label{p:34}
For $\e>0$ small enough and $u=\sum_{i=1}^p\a_i\wtilde{\d_i}_{(a_i,\l_i)}\in V(p,\e)$, the following expansions hold
\begin{align*}
\langle\n J(u),\l_i\frac{\partial \tilde{\d_i}}{\partial \l_i}\rangle_\mathcal{P} = & 2J(u)\left(\frac{n-4}{n} c_2\a_i
\frac{4\D K(a_i)}{\l_i ^2K(a_i)}-c_1\sum_{j\ne i}\a_j\l_i\frac{\partial \e_{ij}}{\partial\l_i}\right)+R\\
\langle\n J(u),\frac{1}{\l_i}\frac{\partial \tilde{\d_i}}{\partial a_i}\rangle_\mathcal{P} = & -2J(u)\left(c_3\a_i\frac{\n K(a_i)}{\l_iK(a_i)}+c_1\sum_{j\ne i}\frac{\a_j}{\l_i}\frac{\partial \e_{ij}}{\partial a_i}\right)+ O\left(\frac{1}{\l_i^2}\right)+ R,
\end{align*}
where $R= o\left( \sum \frac{1}{\l_k ^2}+\sum_{k\ne r}\e_{kr}\right)$.
\end{pro}
\begin{pf}
Using \eqref{ed} and Proposition 2.4 of \cite{BH}, the proof immediately follows from Propositions 3.5 and 3.6 of \cite{BE1}.
\end{pf}
\section{ Characterization of the critical points at infinity }
In this section, we provide the characterization of the
critical points at infinity. First, we construct a special pseudogradient for the associated variational problem for which the Palais-Smale condition is satisfied along the decreasing flow lines, as long as these flow lines do not enter the neighborhood of critical points $y_i$ of $K$ such that $-\D K(y_i)>0$. As a by product of the construction of such a pseudogradient, we are able to determine the critical points at infinity of our problem. 
 \begin{pro}\label{pp:41}
 For $p\geq 2$, there exists a pseudogradient $W$
 so that the following holds.\\
\noindent There is a constant $c>0$ independent of
$u=\sum_{i=1}^p\a_i \wtilde{\d}_i\in V(p,\e)$ so that

$$\langle -\n J(u),W\rangle_{\mathcal{P}} \geq c\biggl(\sum_{i=1}^p\frac{\mid\n
K(a_i)\mid}{\l_i}+\frac{1}{\l_i ^2}+\sum_{i\ne j}\e_{ij}\biggr).\leqno{(a)}
$$
$$\langle-\n J(u+\ov{v}),W+\frac{\partial \ov{v}}{\partial
(\a_i,a_i,\l_i)}(W)\rangle_{\mathcal{P}} \geq c\biggl(\sum_{i=1}^p\frac{\mid\n
K(a_i)\mid}{\l_i}+\frac{1}{\l_i ^2}+\sum_{i\ne j}\e_{ij}\biggr).\leqno{(b)}
$$
$(c)$ $\mid W\mid$ is bounded. Furthermore, $|d \l_i(W)|\leq c\l_i$ for each $i$ and the only case where the maximum of the $\l_i$'s increases along $W$ is when each point $a_i$ is close to a critical point $y_{j_i}$ of $K$ with $-\D K(y_{j_i})>0$ and $j_i\ne j_r$ for $i\ne r$.
\end{pro}
\begin{pf}
We order the $\l_i$'s, for the sake of simplicity we can assume that:
$\l_1\leq ...\leq\l_p$. Let 
$$I_1=\{i|\, \l_i\mid\n K(a_i)\mid \geq C_1'\}, \quad I_2=\{ 1\}\cup \{i\,\mid \l_j \leq M \l_{j-1}, \mbox{ for each } j\leq i\},$$ where $C_1'$ and $M$ are two positive large constants. Set
$$
Z_1=\sum_{i\in I_1}\frac{1}{\l_i}\frac{\partial
\wtilde{\d}_i}{\partial a_i}\frac{\n K(a_i)}{\mid\n K(a_i)\mid}.$$
Using Proposition \ref{p:34}, we derive that
\begin{align}\label{ee:52}
\langle-\n J(u),Z_1\rangle_{\mathcal{P}} &\geq c\sum_{i\in I_1} \frac{\mid\n
K(a_i)\mid}{\l_i}+O\left(\sum _{j\in I_2}\frac{1}{\l_i}\bigg|\frac{\partial
\e_{ij}}{\partial a_i}\bigg|\right)\notag\\
&+O\left(\sum_{i\in I_1}\frac{1}{\l_i ^2}+\sum_{j\notin I_2} \e_{ij}\right)+R.
\end{align}
Observe that, if $j\in I_2$ then 
\begin{eqnarray}\label{ee:53}
  \frac{1}{\l_i}\bigg|\frac{\partial\e_{ij}}{\partial a_i}\bigg|=\l_j|a_i-a_j|\e_{ij}^{(n-2)/(n-4)}=o(\e_{ij}).
\end{eqnarray}
Using also the fact that $i\in I_1$, thus, \eqref{ee:52} becomes
\begin{eqnarray}\label{ee:54}
\langle-\n J(u),Z_1\rangle_{\mathcal{P}} \geq c\sum_{i\in I_1} \frac{\mid\n
K(a_i)\mid}{\l_i}+ \frac{1}{\l_i^2}+O\left(\sum_{j\notin I_2} \e_{ij}\right)+R.
\end{eqnarray}
Now, we will distinguish two cases.\\
{\bf case 1} $I_1\cap I_2\ne \emptyset$.
In this case, we define 
$$Z_2=-M_1\sum_{i\notin I_2}2^i\l_i\frac{\partial
\wtilde{\d}_i}{\partial \l_i}-m_1\sum_{i\in I_2}\l_i\frac{\partial
\wtilde{\d}_i}{\partial \l_i},$$
where $M_1$ is a large constant and $m_1$ is a small constant.\\
Using Proposition \ref{p:34}, we derive 
\begin{align}\label{ee:55}
\langle-\n J(u),Z_2\rangle_{\mathcal{P}} &\geq cM_1\sum_{i\notin I_2}\left( \sum \e_{ij}+O\left(\frac{1}{\l_i^2}\right) +R\right)\notag\\
&+m_1c\sum_{i\in I_2}\left(\sum_{j\in I_2}\e_{ij} +O\left(\frac{1}{\l_i^2}+\sum_{j\notin I_2} \e_{ij}\right) +R\right).
\end{align}
Now, we define $Z_3=Z_1+Z_2$. Using \eqref{ee:54} and \eqref{ee:55}, we derive that
\begin{align}\label{ee:57}
&\langle-\n J(u),Z_3\rangle_{\mathcal{P}}\notag\\
&\quad \geq c\sum_{i\in I_1} \frac{\mid\n
K(a_i)\mid}{\l_i}+ \frac{1}{\l_i^2} + c\sum_{j\ne i} \e_{ij}+O\left( \sum_{i\notin I_2}\frac{M_1}{\l_i^2}+\sum_{i\in I_2}\frac{m_1}{\l_i ^2}\right) +R.
\end{align}
Observe that, since $I_1\cap I_2\ne \emptyset$, we can make $1/\l_k^2$  appear, for $k\in I_2$, in the lower bound of \eqref{ee:57} and therefore all the $\l_i^{-2}$'s can  appear in the lower bound of \eqref{ee:57}. Notice that for $i\notin I_1$, we have $\l_i\mid \n K(a_i)\mid \leq C_1'$. Thus, if we choose $M_1\leq M$ and $m_1 << M^p$,  \eqref{ee:57} becomes
 \begin{eqnarray}\label{ee:58}
\langle-\n J(u),Z_3\rangle_{\mathcal{P}} \geq c\sum_{i=1}^p \frac{\mid\n
K(a_i)\mid}{\l_i}+ \frac{1}{\l_i^2} + c\sum_{j\ne i} \e_{ij}.
\end{eqnarray}
{\bf case 2} $I_1\cap I_2 =\emptyset$.
In this case, for each $i\in I_2$, the point $a_i$ is close to a critical point $y_{k_i}$ of $K$. We claim that $k_i\ne k_j$ for $i\ne j$ that is each neighborhood $B(y,\rho)$, for $\rho$ small enough, contains at most one point $a_i$ with $i\in I_2$. Indeed, arguing by contradiction, let us  suppose  that there exist $i,j\in I_2$ such that $a_i, a_j \in B(y,\rho)$. Since $y$ is nondegenerate we derive that $|\n K(a_k)|\geq c|y-a_k|$ for $k=i,j$ and therefore (we assume that $\l_i\leq \l_j$) $\l_i|a_i-a_j|\leq c$. This implies that $\e_{ij} \geq c (\l_i/\l_j)^{(n-4)/2}$, a contradiction with $\l_i$ and $\l_j$ are of the same order. Thus our claim follows.\\  
Let us introduce 
$$I_3=\{ i\in I_2| \D K(a_i) > 0 \}.$$
{\bf 1st subcase} $I_3 \ne \emptyset$.
In this case we define 
$$Z_4=-\sum_{i\in I_3}\l_i\frac{\partial
\wtilde{\d}_i}{\partial \l_i}-M_1\sum_{i\notin I_2}2^i\l_i\frac{\partial
\wtilde{\d}_i}{\partial \l_i}.$$
Using Proposition \ref{p:34} we derive 
\begin{align}\label{ee:59}
\langle-\n J(u),Z_4\rangle_{\mathcal{P}} & \geq c\sum_{i\in I_3}\left(\frac{1}{\l_i^2}+O\left( \sum \e_{ij}\right)\right)\notag\\
&+M_1c\sum_{i\notin I_2}\left(\sum_{j\ne i}\e_{ij} +O\left(\frac{1}{\l_i^2}\right)\right) +R.
\end{align}
Observe that, if $i,j\in I_2$, we have $|a_i-a_j|\geq c$ then (since $n\geq 9$)
 \begin{eqnarray}\label{ee:60}
\e_{ij}=O\left(\l_i^{-5}+\l_j^{-5}\right).
\end{eqnarray}
For $Z_5= Z_4+Z_1$, using \eqref{ee:54}, \eqref{ee:59}, \eqref{ee:60} and choosing $M_1\leq M$,  we obtain
 \begin{eqnarray}\label{ee:61}
\langle-\n J(u),Z_5\rangle \geq c\sum_{i=1}^p \frac{\mid\n
K(a_i)\mid}{\l_i}+ \frac{1}{\l_i^2} + c\sum_{j\ne i} \e_{ij}.
\end{eqnarray}
{\bf 2nd subcase} $I_3=\emptyset$.
In this case we define 
$$Z_6=\sum_{i\in I_2}\l_i\frac{\partial
\wtilde{\d}_i}{\partial \l_i}-M_1\sum_{i\notin I_2}2^i\l_i\frac{\partial
\wtilde{\d}_i}{\partial \l_i}+Z_1.$$
Using Proposition \ref{p:34}, as in the above subcase, we derive that 
 \begin{eqnarray}\label{ee:62}
\langle-\n J(u),Z_6\rangle \geq c\sum_{i=1}^p \frac{\mid\n
K(a_i)\mid}{\l_i}+ \frac{1}{\l_i^2} + c\sum_{j\ne i} \e_{ij}.
\end{eqnarray}

The vector field $W$ will be a convex combination of all $Z_3$, $Z_5$ and $Z_6$. Thus the proof of claim $(a)$ is completed.\\
By its definition, $W$ is bounded and we have $|d\l_i(W)| \leq c \l_i$ for each $i$. Observe that, the only case when the maximum of the $\l_i$'s increases is where $I_2=\{1,...,p\}$ and $I_1=I_3=\emptyset$, that means each $a_i$ is close to a critical point $y_{j_i}$ of $K$ with $j_i\ne j_r$ for $i\ne r$ and $-\D K(y_{j_i}) > 0$ for each $i$. Hence claim $(c)$ follows. \\
Finally, arguing as in Appendix B of \cite{BCCH}, claim $(b)$ follows from claim $(a)$ and Lemma \ref{l:22}.
    \end{pf}
\begin{pro}\label{p:41}
Let $n\geq 9$. Assume that $J$ has no critical point in $\Sig^+$.
Under the assumptions $(A_1)$ and $(A_2)$, the only
critical points at infinity  under the level $c_\infty{(y_1,y_1)}$ are:
$$
(y_0)_\infty  , \quad (y_1)_\infty \quad\mbox{
 and } \quad (y_0,y_1)_\infty.
$$
Moreover, the Morse indices of such critical points at innfinity are $n- index(K,y_0)$ $=0$, $n-index(K,y_1)$ and $1+n-index(K,y_1)$ respectively.
\end{pro}
\begin{pf}
Using Proposition \ref{p:21}, we derive that $\mid\n J\mid \geq c$ in $\Sig^+ \setminus \cup_{p\geq 1} V(p,\e)$, where $c$ is a positive constant which depends only on $\e$. It only remains to see what happens in $\cup_{p\geq 1} V(p,\e)$. From Proposition \ref{pp:41}, we know that the only region where the maximum of the $\l_i$'s increases along the pseudogradient $W$, defined in Proposition \ref{pp:41}, is the region where each $a_i$ is close to a critical point $y_{j_i}$ of $K$ with $-\D K(y_{j_i})>0$ and $j_i\ne j_r$ for $i\ne r$. In this region, arguing as in \cite{B2}, we can find a change of variables:
$$
(a_1,...a_p, \l_1,...,\l_p) \longrightarrow (\tilde{a}_1,..., \tilde{a}_p,
\tilde\l_1,...,\tilde\l_p):= (\tilde{a},\tilde{\l})
$$
 such that
\begin{align}\label{mokh}
&J\biggl(\sum_{i=1}^p\a_i\tilde\d_{(a_i,\l_i)}+v\biggr)\\
&\quad =\Psi (\wtilde{a},\wtilde{\l}):=
\frac{S_n^{4/n}\sum \a_i^2}{\biggl(\sum\a_i^{\frac{2n}{n-4}}K(\wtilde{a}_i)\biggr)^{\frac{n-4}{n}}}\biggl(1-
(c-\eta)\sum_{i=1}^p \frac{\D K(y_{j_i})}
{\wtilde{\l}_i^2 K(y_{j_i})^{\frac{n}{4}}}\biggr) + \mid V\mid ^2,\notag
\end{align}
where $\eta$ is a small positive constant and $c=c_2 (n-4)/n \left(\sum K(y_{j_i}^{(4-n)/4}\right)^{-1}$, with $c_2$ is defined in Proposition \ref{p:35}.
This yields a split of variables $\tilde{a}$ and $\tilde\l$. Thus it is easy to see that if the $\a_i$'s are in their maximum and  $\wtilde{a}_i=y_{j_i}$ for each $i$, only the $\wtilde{\l}_i$'s can move. To decrease
the functional $J$, we have to increase the $\wtilde{\l}_i$'s, thus we
obtain a critical point at infinity only in this region. It remains to compute the Morse index of such critical points at infinity. For this purpose, we observe that $-\D K(y_{j_i})>0$ for each $i$ and the function $\Psi$ admits on the variables $\a_i$'s an absolute degenerate maximum with one dimensional nullity space and an absolute minimum on the variable $v$. Then the Morse index of such critical point at infinity is equal to $(p-1+ \sum_{i=1}^p (n-index(K,y_{j_i})))$. Thus our result follows.
 \end{pf}

In Proposition \ref{p:41}, we have assumed that $J$ has no critical point in $\Sig^+$. When such an assumption is removed, new critical points at infinity of $J$ appear. Indeed, we have the following result:
\begin{pro}\label{pp:43}
Let $n\geq 9$. Let $w$ be a nondegenerate solution of (1). Then,
$$(y_0,w)_\infty ,\quad (y_1,w)_\infty\quad \mbox{and}\quad (y_0,y_1,w)_\infty$$ are
critical points at infinity. The
Morse index of the critical points are respectively equal to
$$index(w)+1, \,  index(w)+ index((y_1)_\infty)+1 \mbox{ and }
index(w)+index((y_1)_\infty)+2.$$
\end{pro}
The proof of this proposition  immediately follows from
 Proposition \ref{p:41} and the following result:
\begin{pro}\label{c:310}
There is an optimal $(\ov{v},\ov{h})$ and a change of variables
$v-\ov{v}\to V$ and $h-\ov{h}\to H$ such that $J$ reads as
\begin{align*}
J(u)=& \frac{S_n\sum_{i=1}^p\a_i^2+\a_0^2||w||^2}{
(S_n\sum_{i=1}^p\a _i
^{\frac{2n}{n-4}}K(a_i)+\a_0^{\frac{2n}{n-4}}||w||^2
)^{\frac{n-4}{n}}} \left[1-\frac{n-4}{n\b_0}c_2
\sum_{i=1}^p\frac{\a _i ^{\frac{2n}{n-4}}4\D K(a_i)}{\l _i^2}\right.\\
 & \left.-\frac{c_1}{2\gamma_0}\sum_{i\ne j}\a _i\a _j\e _{ij}
 +o\left(\sum_{i\ne j}\e _{ij}
+\sum_{i=1}^p\frac{1} {\l _i^2} \right)\right]
+||V||^2-||H||^2.
\end{align*}
Furthermore, we have the following estimates:
\begin{align*}
||\ov{h}|| & \leq c\sum_i\frac{1}{\l_i
^{(n-4)/2}}\\
||\ov{v}|| & \leq c \biggl[ 
\sum_{i=1}^p\left(\frac{\mid\n K(a_i)\mid}{\l_i}+\frac{1}{\l_i
  ^2}\right) + \sum_{i\ne
  j}\e_{ij}^{\min \left(1, \frac{n+4}{2(n-4)}\right)}(\log\e_{ij}^{-1})^{\min\left(\frac{n-4}{n}, \frac{n+4}{2n}\right)}\biggr].
\end{align*}
\end{pro}
Before giving the proof of Proposition \ref{c:310}, we need to prove the following lemma:
\begin{lem} \label{l:q} The following Claims are true: 
\begin{description}
\item[\it{(a)}]~~$Q_1(v,v)$ is a quadratic form positive definite in\\
 $E_\e '=\{v\in H^2(S^n)|\, v $ $\in T_w(W_s(w)), \mbox{ and }
v \mbox{ satisfies } (W_0)\}$.
\item[\it{(b)}]~~ $Q_2(h,h)$ is a quadratic form negative definite in
$T_w(W_u(w))$.
\end{description}
\end{lem}
\begin{pf}
 Claim (b) follows immediately, since $h\in
T_w(W_u(w))$. Next we are going to prove claim $(a)$.We split $T_w(W_s(w))$ into $E_\g \oplus F_\g$ where $E_\g$ and $F_\g$ are orthogonal for $\langle , \rangle_\mathcal{P}$ and as well as for the quadratic form associated to $w$ and such that
$$
\begin{cases}
\mid\mid v\mid \mid^2 -\frac{n+4}{n-4} \int K w^{8/(n-4)}v^2 \geq (1-\g) \mid\mid v\mid\mid^2\quad\mbox{on}\quad F_\g\\ 
\mbox{dim}( E_\g) < + \infty.
\end{cases}
$$
We choose $\g$ small enough such that $0<\g< \bar\a/4$, where $\bar\a$ is the first eigenvalue of $-\D -\frac{n+4}{n-4}\tilde{\d}_{(a,\l)}^{8/(n-4)}$. Notice that $\bar\a$ is independent on $\tilde{\d}_{(a,\l)}$. Since $dim (E_\g) <\infty$ then we have
$$
\int \tilde\d_i ^{8/(n-4)}v_1^2 =o(||v_1||^2)\quad \forall \,v_1\in E_\g,\,\,\mbox{and } \forall\,i.
$$
 Now let
\be\label{q:1}
v=v_1+v_2,\quad\mbox{with}\quad v_1\in E_\g,\,\,v_2 \in F_\g.
\ee
Then
\begin{align*}
Q_1(v,v)&=||v_1||^2+||v_2||^2 -\sum_{i=1}^p \frac{n+4}{n-4}\int \tilde\d_i^{8/(n-4)}\left(v_1^2+v_2^2 +2v_1v_2\right)\\
&-\frac{n+4}{n-4}\int K w^{8/(n-4)}\left(v_1^2+v_2^2 +2v_1v_2\right)\\
&=||v_1||^2+||v_2||^2 -\sum_{i=1}^p \frac{n+4}{n-4}\int \tilde\d_i^{8/(n-4)}\left(v_1^2+v_2^2\right)\\
&-\frac{n+4}{n-4}\int K w^{8/(n-4)}\left(v_1^2+v_2^2\right)+ o\left(||v_1|| ||v_2||\right)
\end{align*}
This implies that
\begin{align*}
Q_1(v,v)&\geq ||v_1||^2+(1-\g)||v_2||^2 -\sum_{i=1}^p \frac{n+4}{n-4}\int \tilde\d_i^{8/(n-4)}v_2^2 \\
& \qquad -\frac{n+4}{n-4}\int K w^{8/(n-4)}v_1^2+ o\left(||v_1|| ||v_2||+||v_1||^2\right) \\
&\geq (1-\g)||v_2||^2 -\sum_{i=1}^p \frac{n+4}{n-4}\int \tilde\d_i^{8/(n-4)}v_2^2 + o\left( ||v_2||^2\right) + \a'||v_1||^2. 
\end{align*}
It remains to study the term
$$
||v_2||^2 -\sum_{i=1}^p \frac{n+4}{n-4}\int \tilde\d_i^{8/(n-4)}v_2^2.
$$
Observe that $v$ is orthogonal to $span \{\tilde\d_i, \l_i\frac{\partial\tilde\d_i}{\partial \l_i}, \frac{1}{\l_i}\frac{\partial\tilde\d_i}{\partial a_i}, \,\,1\leq i\leq p\}$ but not $v_2$. Since $v_1$ belongs to a finite dimensional space, we have
\be\label{q;2}
\forall \var \in \cup_{i\leq p}\{ \tilde\d_i, \l_i\frac{\partial\tilde\d_i}{\partial \l_i}, \frac{1}{\l_i}\frac{\partial\tilde\d_i}{\partial (a_i)_j}\}, \,\, \mid \langle v_1, \var\rangle_\mathcal{P}\mid \leq ||v_1||_\infty \int |\D^2 \var|=o(||v_1||).
\ee
Now, we write
\be\label{q:3}
v_2=\bar{v}_2 + \sum_{i} A_i \tilde\d_i + \sum_{i} B_i \l_i\frac{\partial\tilde\d_i}{\partial \l_i}+ \sum_{i,j} C_{ij}\frac{1}{\l_i}\frac{\partial\tilde\d_i}{\partial (a_i)_j},
\ee
with  $\bar{v}_2 \in span \{\tilde\d_i, \frac{\partial\tilde\d_i}{\partial \l_i}, \frac{\partial\tilde\d_i}{\partial (a_i)_j}, \,\, i\leq p, \, j\leq n\}^\bot$.\\
Thus, we have (see \cite{BE1})
$$
||\bar{v}_2||^2 -\sum_{i=1}^p \frac{n+4}{n-4}\int \tilde\d_i^{8/(n-4)}\bar{v}_2^2\geq \frac{\bar\a}{2}||\bar{v}_2||^2.
$$
Notice that
\begin{align}\label{q:4}
||v_2||^2 &-\sum_{i=1}^p \frac{n+4}{n-4}\int \tilde\d_i^{8/(n-4)}v_2^2= ||\bar{v}_2||^2 + O\left(\sum_i A_i^2 + B_i^2 + \sum_j C_{ij}^2\right)\notag\\
&-  \sum_{i=1}^p \frac{n+4}{n-4}\int \tilde\d_i^{8/(n-4)}\bar{v}_2^2 + O\left(||\bar{v}_2||(|A_i| + |B_i| + \sum_j |C_{ij}|)\right)
\end{align}
Using \eqref{q:1}-\eqref{q:3}, we obtain
$$
A_i=o(||v_1||),\quad B_i=o(||v_1||)\quad\mbox{and } C_{ij}=o(||v_1||)\,\mbox{ for each } i,j.
$$
Thus, using \eqref{q:4}, we derive that
$$
Q_1(v,v)\geq -\g||v_2||^2 + \frac{\bar\a}{2}||\bar{v}_2||^2 + o\left(||v_1||^2 + ||v_2||^2\right) + \a'||v_1||^2.
$$
But
$$
||v_2||^2= ||\bar{v}_2||^2 +  O\left(\sum_i A_i^2 + B_i^2 + \sum_j C_{ij}^2\right)= ||\bar{v}_2||^2 + o(||v_1||^2).
$$
Thus
$$
Q_1(v,v)\geq \left(\frac{\bar\a}{2}-\g\right)||v_2||^2 + \a'||v_1||^2 + o\left(||v_1||^2 + ||v_2||^2\right).
$$
Since $\g < \bar\a/4$, claim $(a)$ follows. The proof of our lemma is thereby
completed.\end{pf}\\
\begin{pfn}{Proposition \ref{c:310}}
 By Proposition \ref{p:35} the expansion of $J$ with respect to $h$
(respectively to $v$) is very close, up to a multiplicative
constant, to $Q_2(h,h)+f_2(h)$ (respectively $Q_1(v,v)-f_1(v)$).  By Lemma \ref{l:q} there is a unique maximum $\ov{h}$ in the space of $h$
(respectively a unique minimum $\ov{v}$ in the space of $v$).
Furthermore, it is easy to derive that $||\ov{h}||\leq
c||f_2||=O(\sum_i \l_i ^{(4-n)/2})$ and $||\bar{v}||\leq c||f_1||$. The estimate of $\ov{v}$
follows from Lemma \ref{l:22}. Then our result
follows.
\end{pfn}
\section{ Proof of Theorems}
Let us start by proving the
following results.
\begin{pro}\label{p:61}
Let $z_1,z_2\in X$ be such that $-\D K(z_i)>0$
for $i=1,2$, $z_1\ne z_2$ and $z_1$, $z_2$ satisfy assumption $(A_2)$. If we assume\\  
(a) $\quad
J(\frac{1}{K(z_1)^{(n-4)/8}}\tilde\d_{(z_1,\l)}+\frac{1}{K(z_2)^{(n-4)/8}}
\tilde\d_{(z_2,\l)})\geq c_{\infty}(z_1,z_2)+\d$,\\
(b) $\quad (\partial/\partial\mu)
J(\frac{1}{K(z_1)^{(n-4)/8}}\tilde\d_{(z_1,\mu)}+\frac{1}{K(z_2)^{(n-4)/8}}
\tilde \d_{(z_2,\mu)})\big|_{\mu=\l}<0$,\\ 
then $I(V)=0$ for any $V\in \mathcal{F}$.
\end{pro}
\begin{pf}
An abstract topological argument displayed in \cite{B2}, pages 358--369, which extends to our framework,  shows that the value of $I(V)$ is constant for any $V \in \mathcal{F}$. Now, let $\e>0$ and $K_\e=1+\e K$. Let $J_\e$ be the associated variational problem. As $\e$ tends to zero, $J_\e$ tends to $J_0$ in the $C^1$ sense, where $J_0$ is the functional defined  replacing $K$ by $1$ in \eqref{functional}. On the other hand, using Proposition \ref{p:35}, we see that 
$$
J_\e (\a_1 \tilde\d_{(a_1,\l)}+\a_2 \tilde\d_{(a_2,\l)})\leq 2S^{4/n}\left(1-\frac{c}{\l^{n-4}}+O(\e)\right),
$$
where $c$ is independent of $\e$ and $2S^{4/n}$ is the level to which a critical point at infinity of $2$ masses of $K_\e$ converges when $\e\to 0$. Thus, we can assume  $\e$ is so small that all critical points at infinity of $J_\e$ (of two masses or more) are above $f_\l(B_2(X))$. Clearly, for $\e$ small, $C_\d(z_1,z_2)$ is above $(2S^{4/n} + \d/2)$. We derive that
$$
W_u^\e(f_\l(B_2(X))) . C_\d(z_1.z_2)=0.
$$
Notice that, decreasing $\l$, we complete a homotopy of  $f_\l(B_2(X))$ that increases the interaction of any masses, and therefore remains below $C_\d(z_1,z_2)$. This implies that for each $\mu\in [1,\l]$ we have
$$
W_u^\e(f_\mu(B_2(X))) . C_\d(z_1.z_2)=0.
$$
Recall that
\be\label{invariant}
I(V)=\tau+\sum_{w_{2k+1}\in
A_{2k+1}}(f_\l(B_2(X)).w_{2k+1})((y_0,w_{2k+1})_\infty .C_\d).
\ee
Thus, we need to compute $f_\l(B_2(X)).w_{2k+1}$ for any $w_{2k+1}\in A_{2k+1}$. Let
$$
F=\cup_{\mu =1}^\l f_\mu (B_2(X)).
$$
We can assume that $F$ is a compact manifold in dimension $2k+2$. The singularity of $F$ is $\cup_{\mu=1}^\l f_\mu(B_1(X))$ which is of a dimension less than $(k+1)$, this singularity cannot dominate $w_{2k+1}$. We deduce that $F\cap \bar{W}_s(w_{2k+1})$ is a compact manifold of dimension one. Thus the cardinal of $\partial (F\cap \bar{W}_s(w_{2k+1}))$ is equal to zero, where $\partial$ is the boundary homomorphism of $S_{2k+2}(\Sig^+)$.\\
Observe that 
$$
\partial F= f_1(B_2(X)) + f_\l(B_2(X)).
$$
It follows that
$$
f_\l(B_2(X)).w_{2k+1} = f_1(B_2(X)).w_{2k+1} + F. \partial^{-1}\left(W_s(w_{2k+1})\right).
$$
Along this homtopy, the trace of $f_\mu(B_2(X))$ might intersect, for some values, $\partial^{-1}\left(W_s(w_{2k+1})\right)$, where $\partial^{-1}\left(W_s(w_{2k+1})\right)$ is made of stable manifolds of critical points of index $2k+2$. Therefore the abstract argument of \cite{B2} applies, and the invariant remains unchanged. For $\mu =1$ at the end of the homotopy $B_2(X)$ is mapped onto a single function and $\left(f_1(B_2(X)). w_{2k+1}\right)$ is therefore zero. Thus, $I(V)$ at the end of the homotopy is equal to zero, and the results follow.
\end{pf}

Now, we are going to prove Theorem \ref{t:1}.

\begin{pfn} {Theorem \ref{t:1}}
Arguing by contradiction, we assume that $J$ has no critical point in $\Sig^+$.
It follows from Proposition \ref{p:41} that $A_{2k+1} =\emptyset$. Therefore combining \eqref{invariant}, Proposition \ref{p:61} and the fact that $\tau =1$, we derive a contradiction. The proof of our result is thereby completed.
\end{pfn}

The sequel of this section is devoted to the proof of Theorem \ref{t:2}.

\begin{pfn} {Theorem \ref{t:2}}
 In the sequel, we denote by
$\Pi _a$ the  stereographic projection through a point $a \in S^n$. This projection induces an isometry $i : H^2(S^n) \to
\mathcal{H}(\R^n) $ according to the following formula
$$
(iv)(x)= \left(\frac{2}{1+|x|^2}\right)^{(n-4)/2}v(\Pi
_a^{-1}(x)), \qquad v\in H^2(S^n), \, x\in \R^n,
$$
where $\mathcal{H}= \{u\mid u\in L^{2n/(n-4)}(\R^n),\,\, \D u\in L^2(\R^n)\}$.
Now, let $a$ in $ S^n$ (it is easy to see that
$\pi_{-a}(a)=o$ and $i(\tilde\d_{(a,\l)})=\d_{(o,\l)}$).

Let $a_1$, $a_2$ in $ S^n$ and $\rho_1$, $\rho_2$ be two
positive constants (we choose $\rho_1$ and $\rho_2$ such that
$B(a_1,\rho_1')\cap B(a_2,\rho_2')$ is empty i.e.
$\rho_1'+\rho_2'< d(a_1,a_2)$). Let
$$
u=\a_1\tilde\d_{(a_1,\l_1)}+
\a_2\tilde\d_{(a_2,\l_2)}+v,\quad \mbox{  with } \quad \a_i=K(a_i)^{(4-n)/8}$$
where $v$ satisfies $(V_0)$ which is defined in \eqref{V0}.

 We now write down the expansion of  $J(u)=N/D$
with
\be\label{N'}
N =  S_n\sum_{i=1,2}\frac{1}{K(a_i)^{(n-4)/4}}+||v||^2 
 +O\left(\sum_{i=1,2}
\frac{1}{K(a_i)^{(n-4)/4}}\frac{1}{(\l_i\rho_i)^{n-4}}\right),
\ee
\begin{align}\label{D}
D^{\frac{n}{n-4}} 
&=\sum_{i=1}^2\frac{1}{K(a_i)^{\frac{n}{4}}}\int K\tilde\d_i ^{\frac{2n}{n-4}}+\frac{2n}{n-4}\int K \left(\sum\a_i\tilde\d_i\right)^{\frac{n+4}{n-4}}v\notag\\
& + \frac{n(n+4)}{(n-4)^2}\int K\left(\sum \a_i
\tilde\d_i\right) ^{\frac{8}{n-4}}v^2+\sum O\left(\frac{1+ R_{K,i}^2 }{K(a_i)^{(n-4)/4}(\l_i\rho_i)^{n-4}}\right)  \notag\\
& +O\left(\sup_{S^n}K\bigg(||v||^{\frac{2n}{n-4}}+(if\,
n < 12)\frac{||v||^3}{K(a_i)^{\frac{12-n}{8}}}\bigg)\right).
\end{align}
where $R_{K,i}$ satisfies
\be \label{R}
R_{K,i}=\frac{|\n K(a_i)|}{\l_iK(a_i)}+\frac{|D^2 K(a_i)|}{\l_i^2K(a_i)}+\sup_{B_i}\frac{|D^3 K|}{\l_i ^3K(a_i)}.
\ee
 Now, assuming $\l_i$ and $\l_i\rho_i$ are large, we write
\begin{align*}
\int_{S^n}  K \tilde\d_i ^{\frac{2n}{n-4}} = &
 K(a_i) S_n + \frac{4\D K(a_i)}{\l_i^2}\left(c_2+O\left(\frac{1}{(\l_i\rho_i)^{n-2}}\right)\right)\\
 & +O\left(\sup_{B_i}\frac{|D^3 K|}{\l_i^3}+ \frac{\sup K}{(\l_i\rho_i)^{n}}\right).
\end{align*}

Thus
\begin{align}\label{5.43}
J(u) =&  \biggl(S_n\sum_{i=1}^2\frac{1}{K(a_i)^{\frac{n-4}
{4}}}\biggr)^{4/n} \bigg[1-\frac{c_2(n-4)}{n\beta}\sum_{i=1}^2\frac{4\D K(a_i)}{\l_i^2 K(a_i)^{n/4}}\\
 & +O\left(\frac{1}{\b}\sum_{i=1}^2\frac{1+R_{K,i}^2}{K(a_i)^{\frac{n-4}{4}}
(\l_i\rho_i)^{n-4}}\right)  -\frac{1}{\b}f(v)\notag   \\
& +\frac{1}{\b}\biggl(||v||^2-\frac{n+4}{n-4}\int K\left(\sum \a_i
\tilde\d_i\right)^{\frac{8}{n-4}}v^2\biggr)\notag\\
 &  +\sum \frac{1}{\b K(a_i)^{n/4}}O\biggl(\frac{|\D K(a_i)|}{\l_i^2(\l_i\rho_i)^{n-2}}+\sup_{B_i}\frac{|D^3K|}{\l_i ^3}+ \sup_{S^n} K \frac{1}{(\l_i \rho_i)^n} \biggr)\notag\\
 & +\frac{1}{\b}\sum_{i=1}^2
O\biggl(\sup_{S^n}K \biggl(||v||^{\frac{2n}{n-4}}+(if\,
 n<12)\frac{||v||^3}{K(a_i)^{(12-n)/8}}\biggr)\biggr)\bigg],\notag
\end{align}
where $\b=S_n\sum_{i=1}^2{1/K(a_i)^{(n-4)/4}}$ and where
$$
f(v)=2 \int_{S^n} K (\a_1\tilde\d_1+\a_2\tilde\d_2)^{\frac{n+4}{n-4}}v.
$$
Notice that 
\begin{align}\label{540}
f(v) & =2\sum \a_i^{\frac{n+4}{n-4}}\int_{S^n}K \tilde\d_i ^{\frac{n+4}{n-4}}v
 + O\left(\int K \sup {}^{\frac{8}{n-4}}(\a_1\tilde\d_1,\a_2\tilde\d_2)\inf(\a_1\tilde\d_1,\a_2\tilde\d_2)|v|\right)\notag  \\
 & = O\left(||v||\sum \frac{1}{K(a_i)^{\frac{n-4}{8}}}\left(R_{K,i}+\frac{\sup K}{K(a_i)} \frac{\log(\l_i\rho_i)^{(n+4)/n}}{(\l_i\rho_i)^\frac{n+4}{2}}\right)\right)
\end{align}
On the other hand, we know from Proposition 3.4 of \cite{BE1} that the quadratic form
\begin{eqnarray}
||v||^2-\frac{n+4}{n-4}\sum_{i=1}^2\int_{S^n}
\tilde\d_i ^{\frac{8}{n-4}}v^2
\end{eqnarray}
is bounded below by $\a_0||v||^2$, $\a_0$ is a fixed constant, on all
$v$'s satisfying $(V_0)$.
 Observe now 
\begin{align} 
\int K & \left(\sum \a_i \tilde \d_i\right)^{\frac{8}{n-4}}v^2  =\sum \int \frac{K}{K(a_i)}\tilde\d_i ^\frac{8}{n-4}v^2 +O\left(\int K (\a_1\tilde\d_1\a_2\tilde\d_2)^\frac{4}{n-4}v^2\right)\notag\\
 &\quad =\sum \int \tilde\d_i^{8/(n-4)}v^2 +O\left(||v||^2\left(\sum \frac{\sup K}{K(a_i)}\frac{\log^{8/n}(\l_i\rho_i)}{(\l_i\rho_i)^4}+R_{K,i}\right)\right)
\end{align}
Thus, if we assume that 
\be\label{546}
\sum \frac{\sup K}{K(a_i)}\frac{\log^{8/n}(\l_i\rho_i)}{(\l_i\rho_i)^4}+R_{K,i}
\ee
is small,  then the quadratic form which comes out of the expansion
\begin{eqnarray}
||v||^2-\frac{n+4}{n-4} \int
K(\a_1\tilde\d_1+\a_2\tilde\d_2)^{\frac{8}{n-4}}v^2
\end{eqnarray}
is definite positive, bounded below by $(\a_0/4)||v||^2$ for $v$
satisfying $(V_0)$. Therefore the functional
\begin{eqnarray}
-f(v)+||v||^2-\frac{n+4}{n-4}\int K\left( \a_1\tilde\d_1+\a_2\tilde\d_2\right)^{\frac{8}{n-4}}v^2
\end{eqnarray}
has a unique minimum $\tilde{v}$ and we have
$||\tilde{v}||=O(||f||)$.\\
 The function $J(u)$ has in fact one more term depending on $v$ which is
\be
\sum_{i=1}^2 O\biggl(\sup_{S^n}K  \biggl(||v||^{
\frac{2n}{n-4}}+(if\, n<12)\frac{||v||^3}{K(a_i)^{(12-n)/8}}
 \biggr)\biggr).
\ee
$J$ is twice differentiable. Therefore, this remainder term is
also twice differentiable and its second differential is easily
checked to be
\be
\sup_{S^n} K O(||v||^{8/(n-4)}) + \sum \frac{\sup K}{K(a_i)^{(12-n)/8}}O(||v||)(\mbox{if } n < 12).
\ee
Thus, if we assume that $(\sup K)O(||f||^{{8}/(n-4)})
\leq \tilde{\tilde{c}}$ (for $n\geq 12$) and (for $n< 12$)\\
$\sup K(K(a_1)^{(n-12)/8}+K(a_2)^{(n-12)/8})O(||f||)\leq \tilde{\tilde{c}}$ where $\tilde{\tilde{c}}$ is a small constant, the functional
\begin{align*}
-f(v)+ & ||v||^2-  \frac{n+4}{n-4}\int  K(x)(\sum \a_i\tilde\d_i)^{\frac{8}{n-4}}v^2 +(\sup K)O(||v||^{\frac{2n}{n-4}})\\
&+(if\, n<12)\sup K (K(a_1)^{(n-12)/8}+K(a_2)^{(n-12)/8})O(||v||^3)
\end{align*}
will have a unique minimum $\ov{v}$ near the origin and it satisfies also
$||\ov{v}||=O(||f||)$. Let us introduce the following neighborhood $V$ of functions $v\in H^2(S^n)$ such that $v$ satisfies $(V_0)$ and  
\be\label{H}
\begin{cases} 
 ||v-\ov{v}||<\frac{\tilde
{c_1}}{(\sup K)^{(n-4)/8}}\quad (\mbox{if } n\geq 12)\\
 ||v-\ov{v}||<\frac{\tilde{c_1}}{\sup K (K(a_1)^{(n-12)/8}+K(a_2)^{(n-12)/8})}\quad (\mbox{ if }n<12).
\end{cases}
\ee
Requiring $v$ to belong to $V$, we let by $\ov{u}=\sum
(1/K(a_i)^{(n-4)/8})\tilde\d_i+\ov{v}$. Then 
\begin{eqnarray}
J(u)=J(\ov{u})+\left(S_n\sum_{i=1}^2\frac{1}{K(a_i)^{(n-4)/4}}
\right)^{(4-n)/n} Q(v-\ov{v},v-\ov{v}),
\end{eqnarray}
where $Q$ is a definite positive form, bounded below by
$(\a_0/4)||v-\ov{v}||^2$ on $V$. An expansion of $J(\ov{u})$ is
easily derived by setting $v=\ov{v}$ in the expansion of $J(u)$ (see \eqref{5.43}) and
using the estimate of $\ov{v}$. Thus,
\begin{align*}
J(\ov{u}) &
=\b^{4/n}\biggl[1-\frac{c_2(n-4)}{n\beta}\sum_{i=1}^2\frac{4\D K(a_i)}{\l_i^2
K(a_i)^{n/4}}+O(\frac{1}{\b}||f||^2)\\
 &
 +O\left(\sum_{i=1}^2\frac{1+R_{K,i}^2}{\b K(a_i)^{\frac{n-4}{4}}
(\l_i\rho_i)^{n-4}}\right)\\
 & +\sum_{i=1}^2\frac{1}{\b K(a_i)^{n/4}}O\biggl(
\frac{\sup K}{(\l_i\rho_i)^{n}} +\frac{|\D
K(a_i)|}{\l_i^2(\l_i\rho_i)^{n-2}} +\sup_{B_i}\frac{|D^3K|}{\l_i
^3}\biggr)\biggr].
\end{align*}
As in Proposition \ref{p:34} and in Appendix B of \cite{BCCH}, we obtain
\begin{align}
\l_j\frac{\partial J(\ov{u})}{\partial\l_j} &
=\b^{\frac{4-n}{n}}\biggl[\frac{8c_2(n-4) \D K(a_j)}{n\l_j^2K(a_j)^{n/4}}
+O\biggl(\sum_{i=1}^2\frac{1}{K(a_i)^{\frac{n-4}{4}}}
\left(\frac{1+R_{K,i}^2}{(\l_i\rho_i)^{n-4}}\right.\\
& \left.+\frac{\sup K}{K(a_i)}\frac{1}{(\l_i \rho_i)^{n}}+\frac{|\D K(a_i)|}{\l_i^2K(a_i)(\l_i\rho_i)^{n-2}} +\sup_{B_i}\frac{|D^3K|}{\l_i^3 K(a_i)}\right)+||f||^2 \biggr)\biggr].\notag
\end{align}
Thus for $\b_1$, $\b_2 \geq 0$, $\b_1+\b_2=1$ and using the estimate of $||f||$ (see \eqref{540}), we derive 
\begin{align}
&\sum_{j=1}^2 \b_j \l_j\frac{\partial J(\ov{u})}{\partial\l_j}
=\b^{\frac{4-n}{n}}\biggl[\frac{8c_2(n-4)}{n} \sum_{j=1}^2\frac{\b_j \D K(a_j)}{\l_j^2K(a_j)^{n/4}}\\
&\quad+\sum_{j=1}^2\b_j O\biggl(\frac{1}{K(a_j)^{\frac{n-4}{4}}}\biggl(
\frac{1+R_{K,j}^2}{(\l_j\rho_j)^{n-4}}+\frac{\sup K}{K(a_j)}\frac{1}{(\l_j\rho_j)^n} +\frac{|\D K(a_j)|}{\l_j^2K(a_j)(\l_j\rho_j)^{n-2}}\notag\\
&\quad +\sup_{B_j}\frac{|D^3K|}{\l_j^3K(a_j)}+R_{K,i}^2+\frac{\sup K^2}{K(a_j)^2}\frac{\log(\l_j\rho_j)^{2(n+4)/n}}{(\l_j\rho_j)^n+4}
\biggr)\biggr) \biggr].\notag
\end{align}
This derivative will remain negative as long as, for a suitable
universal constant $c'_1$, we have for $i=1,2$
\begin{align}\label{I}
 \frac{1}{(\l_i\rho_i)^{n-4}} & +\frac{\sup K}{K(a_i)}\frac{1}{(\l_i\rho_i)^n}+\frac{\sup K^2}{K(a_i)^2}\frac{\log(\l_i\rho_i)^{2(n+4)/n}}{(\l_i\rho_i)^{n+4}}\notag \\ 
 & + \frac{|\n K(a_i)|^2}{\l_i
^2K(a_i)^2}+\sup_{B_i}\frac{|D^3 K|}{\l_i ^3K(a_i)}+ \frac{|D^2 K(a_i)|^2}{\l_i^4 K(a_i)^2} \leq c'_1\frac{-\D K(a_i)}{\l_i^2K(a_i)}.
\end{align}
Taking $c'_1$ to be smaller, if necessary, we derive that, under \eqref{I}
and if $v\in V$, $J(u)$ is bounded below as follows:
\begin{eqnarray}\label{643}
J(u)\geq \b^{4/n}\biggl[1+\frac{1}{\b}\biggl(c_2'\sum_{i=1}^2
\frac{-\D K(a_i)}{\l_i^2
K(a_i)^{n/4}}+\frac{\a_0}{4}||v-\ov{v}||^2\biggr)\biggr].
\end{eqnarray}
To \eqref{I}, other conditions which we used earlier are to be added,
namely
\begin{align}
& ||f||\sup K\left(\sum K(a_i)^{(n-12)/8}\right) \leq c_1''\, \mbox{ if } n<12\label{J1} \\
& ||f|| \sup K^{(n-4)/8} \leq c_1'' \, \mbox{ if } n\leq 12\label{J2}\\
& \frac{\sup K}{K(a_i)}\frac{\log(\l_i\rho_i)^{(n+4)/8}}{(\l_i\rho_i)^4}+R_{K,i} \leq c_1'' \, \mbox{ for } i=1,2.\label{J3}
\end{align}
Finally, all the quantities involved in \eqref{I}, up to the factor
$1/\b$, should be small for the expansions to hold, which amounts
to
\be\label{K}
\frac{1}{\b}\biggl(\sum_{i=1}^2\frac{-\D K(a_i)}{\l_i^2K(a_i)^{n/4}}\biggr)<c'_1.
\ee
We will take
\be\label{A}
 a_1\in \nu^+(z_1),\, \, \nu^+(z_1) \mbox{ small enough so that }
K(z_1)\leq K(a_1)\leq 2 K(y_1).
\ee
We will ask that
\be\label{B}
 a_2\in \nu^+(z_2),\, \, \nu^+(z_2)~ be \mbox{ small enough so that }
K(z_2)\leq K(a_2)\leq 2K(y_1).
\ee
and that
\be\label{C}
\frac{1}{2}K(y_1)\leq K(z_1), \quad \frac{1}{2}K(y_1)\leq K(z_2).
\ee
From \eqref{I} and \eqref{K}, it is easy to derive that $R_{K,i}$ is small. Observe also that, using \eqref{I} and \eqref{J3}, \eqref{I} can be simplified. Finally,  \eqref{I}, \eqref{J1}--\eqref{K} therefore reduce to (after reducing $c'_1$)
\be\label{L}
\begin{cases}
  \frac{1}{\rho_i^2(\l_i\rho_i)^{n-6}}+\frac{|\n K(a_i)|^2}{
K(a_i)^2}+ \frac{|D^2K(a_i)|^2}{\l_i^2 K(a_i)^2} +\sup_{B_i}\frac{|D^3 K|}{\l_i K(a_i)}\leq c'_1\frac{-\D K(a_i)}{K(a_i)}\\
(\sup K/ K(y_1))^{\max(1,(n-4)/8)}||f||K(y_1)^{(n-4)/8} \leq c_1''\\
  \sum_{i=1}^2\frac{|\D K(a_i)|}{\l_i^2 K(a_i)} \leq c'_1\\
\sum (\sup K) K(y_1)^{-1}\log(\l_i\rho_i)^{(n+4)/n}(\l_i\rho_i)^{-4}\leq c_1''.
\end{cases}
\ee
The third condition of \eqref{L} follows from the first one, since $|D^2
K(a_i)|$ dominates $|\D K(a_i)|$ (up to a modification of $c'_1$). Thus
\be\label{L'}
\begin{cases}
  \frac{1}{\rho_i^2(\l_i\rho_i)^{n-6}}+\frac{|\n K(a_i)|^2}{
K(a_i)^2}+ \frac{|D^2K(a_i)|^2}{\l_i^2 K(a_i)^2} +\sup_{B_i}\frac{|D^3 K|}{\l_i K(a_i)}\leq c'_1\frac{-\D K(a_i)}{K(a_i)}\\
(\sup K/ K(y_1))^{\max(1,(n-4)/8)}||f||K(y_1)^{(n-4)/8} \leq c_1''\\
\sum (\sup K) K(y_1)^{-1}\log(\l_i\rho_i)^{(n+4)/n}(\l_i\rho_i)^{-4}\leq c_1''.
\end{cases}
\ee
At this point, following the proof of \cite{B2}, we explain how we
will proceed with the proof of Theorem \ref{t:2}. We wish to compute
$W_u(f_\l(B_2(X))).C_\d(z_1,z_2)$.

Let us define
\begin{eqnarray}
g_\l: B_2(X)\to \Sig ^+,\quad (\a_1,\a_2,a_1,a_2)\to
\frac{\a_1\tilde\d_{(a_1,\l)}+\a_2\tilde\d_{(a_2,\l)}+\ov{v}}
{||\a_1\tilde\d_{(a_1,\l)}+\a_2\tilde\d_{(a_2,\l)}+\ov{v}||}.
\end{eqnarray}
$g_\l$ and $f_\l$ are homotopic (see \cite{B2}). Using also the
fact that $-\D K(z_1)$ and $-\D  K(z_2)$ are
positive, we can choose $\d$ so small such that
\begin{eqnarray}
g_\l( B_2(X)).W_s(C_\d(z_1,z_2))=f_\l( B_2(X)).W_s(C_\d(z_1,z_2)).
\end{eqnarray}
We can accordingly modify $C_\d(z_1,z_2)$ as follows:
\begin{eqnarray}
\tilde C_\d(z_1,z_2)=\tilde\Gamma_{\e_1}(z_1,z_2)\cap
J^{-1}(c_\infty(z_1,z_2)+\d),
\end{eqnarray}
where
\begin{align}
\tilde\Gamma_{\e_1}&(z_1,z_2)=\bigg\{ 
\sum_{i=1,2}\frac{\tilde\d_{(z_i+h_i, \l_i)}}
{K(z_i+h_i)^{\frac{n-4}{8}}}+v/v\in H^2(S^n) \mbox{
satisfies }(V_0),\notag\\
 &|| v-\ov{v}||<\e_1,\,\,
  \l_i>\e_1^{-1}\mbox{ for }i=1,2,\, h_i\in
\nu^+(z_i),\, \mid h_1\mid^2+\mid h_2\mid^2<\e_1\bigg\}.\notag
\end{align}
Clearly, $C_\d(z_1,z_2)$ and $\tilde{C}_\d(z_1,z_2)$ can be
deformed, one into another, using an isotopy above the level
$c_\infty(z_1,z_2)$. Thus
\begin{eqnarray}
g_\l( B_2(X)).W_s(C_\d(z_1,z_2))=\tau(z_1,z_2)=
f_\l(B_2(X)).W_s(C_\d(z_1,z_2)).
\end{eqnarray}
Computing $\tau(z_1,z_2)$ now becomes a matter of defining a
pseudogradient such that the Palais-Smale condition ((P.S.) for
short) is satisfied along decreasing flow lines away from the
critical points at infinity and computing $\tau(z_1,z_2)$ for this
flow. In the absence of solutions, $\tau$ does not depend on this
pseudogradient as long as the asymptotes are as expected. We can
therefore compute $\tau$ with a special flow worrying only about the fact
that it belongs to $\mathcal{F}$ and is admissible. Observe now that, if
we take $\d$ very small, $h_1$ and $h_2$ are as small as we may
wish in $\tilde{C}_\d(z_1,z_2)$ ($\e_1$ has been chosen very
small before $\d$, $\d$ is then chosen so small that
$\tilde{C}_\d(z_1,z_2)$ is a Fredholm manifold of codimension
$2k+2$).\\
To construct the vector field, we need that $(\l_1,\l_2)\in
[A_1,+\infty)\times [A_2,+\infty)$,\\
$(a_1,a_2)\in \nu^+(z_1)\times\nu^+(z_2)$ such that (see \eqref{H} for the definition of  $V$):
\begin{enumerate}
 \item $ B(a_1,\rho_1)\cap B(a_2,\rho_2)=\emptyset$ for each
$(a_1,a_2)\in \nu^+(z_1)\times \nu^+(z_2)$ such that
$c_\infty(a_1,a_2)\leq c_\infty(y_1,y_1)$.
 \item on $\partial([A_1,+\infty)\times [A_2,+\infty)\times V)$,
$$
J(\tilde \d_{(a_1,\l_1)}/K(a_1)^{\frac{n-4}{8}}+\tilde
\d_{(a_2,\l_2)} /K(a_2)^{\frac{n-4}{8}}+v)\geq c_\infty(y_1,y_1),
$$
for any $(a_1,a_2)\in \nu^+(z_1) \times\nu^+(z_2)$.
 \item  \eqref{L'} is satisfied on $[A_1,+\infty)\times
[A_2,+\infty)\times\nu^+(z_1) \times\nu^+(z_2).$
\end{enumerate}

Assuming now that 1), 2) and 3) hold and taking $\l\geq
\max(A_1,A_2)$, we first observe that the expansion of $J$ splits
completely the variable $(\l_1,\l_2)$ from $v-\ov{v}$. Therefore,
we can build our pseudogradient independently on both variables.
In the $(v-\ov{v})$-space, we simply increase $v-\ov{v}$
directionally, if it is non zero, that is
\begin{eqnarray}
\frac{\partial}{\partial s}(v-\ov{v})=v-\ov{v}.
\end{eqnarray}
This increasing component of the pseudogradient will not move the
concentration and will bring the $v$'s on $\partial  V$, if
$v-\ov{v}$ is non zero initially, hence above $c_\infty(y_1,y_1)$. Since
$g_\l(B_2(X))$ is below $c_\infty(y_1,y_1)$, $\tilde C_\d(z_1,z_2)$ and
$g_\l(B_2(X))$ will not intersect through these flow lines. Thus,
any intersection will come from $v=\ov{v}$.

In the case where $c_\infty(a_1,a_2)\leq c_\infty(z_1,z_2)+\d/4$,
in the $(\l_1,\l_2)$-space when $v=\ov{v}$, an increasing pseudogradient can be obtained by decreasing both $\l_1$ and $\l_2$ and
keeping the ratio $\l_1/\l_2$ unchanged (using  condition \eqref{L'}).
The Palais-Smale condition will be satisfied on the decreasing flow lines of such pseudogradient which is defined as such only above $\tilde
C_\d(z_1,z_2)$ and has to be extended to the other regions because,
if any of $\l_1$ or $\l_2$ tends to $+\infty$, then, since the
ratio is unchanged, both tend to $+\infty$ and $J$ (since
$v=\ov{v}$) tends to $c_\infty(a_1,a_2)$ which is below
$c_\infty(z_1,z_2)+\d/4$. However, under the level
$c_\infty(z_1,z_2)+\d/2$ we can construct our pseudogradient such
as we did in Proposition \ref{pp:41}. This one will satisfy the Palais-Smale condition on
decreasing flow lines away from the critical points at infinity
announced in Proposition \ref{pp:41}. Thus, with this suitable extension,
we can freely define, above $c_\infty(z_1,z_2)+\d$, our pseudogradient by decreasing $\l_1$ and $\l_2$ and by taking the ratio
unchanged.

In the other case, which is $c_\infty(a_1,a_2)\geq
c_\infty(z_1,z_2)+\d/4$, this forces $(a_1,a_2)$ in
$\nu^+(z_1)\times \nu^+(z_2)$ to be away from $(z_1,z_2)$, sizeably
away. We can then move $(a_1,a_2)$ in the outwards direction in
$\nu^+(z_1)\times \nu^+(z_2)$. $c_\infty(a_1,a_2)$ then increases,
until it reaches the level $c_\infty(y_1,y_1)$. Since $\l_1$ and $\l_2$ can
be assumed as large as we may wish, this builds a pseudogradient
for $J$ between the level of $C_\d(z_1,z_2)$ and $c_\infty(y_1,y_1)$, in
the region where $\l_1$ and $\l_2$ are extremely large, which
satisfies (P.S.) since the concentration remains unchanged.
 Clearly, we will intersect $g_\l(B_2(X))$ only once, when
$\l_1=\l_2=\l$. The intersection of $g_\l(B_2(X))$ and $W_s(\tilde
C_\d(z_1,z_2))$ then becomes transversal. 

We now need to prove that we can find $A_1$ and $A_2$ such that 2)
holds. 
Assuming that
\begin{eqnarray}
\min \biggl(\frac{K(y_1)}{K(a_1)},\frac{K(y_1)}{K(a_2)}\biggr)\geq
1-c'_0,
\end{eqnarray}
$c'_0$ being a small fixed constant, we can modify the lower-bound in \ref{643} as
follows
\begin{align}\label{m}
J(u) \geq & c_\infty(y_1,y_1)
\biggl(1+ c \biggl(1-\frac{K(a_1)+ K(a_2)}{2K(y_1)}\\
&+\sum_{i=1}^2
\frac{-\D K(a_i)}{\l_i^2K(a_i)}+
\frac{\a_0}{4}K(y_1)^{\frac{n-4}{4}}||v-\ov{v}||^2\biggr)\biggr).\notag
\end{align}
Under \eqref{m}, the set $V$ in \eqref{H} can be replaced by
\begin{eqnarray}
\tilde V =\{ v/\, \, ({\a_0}/{4})K(y_1)^{(n-4)/4}||v- \ov{v}||^2
\leq \tilde c_2\}.
\end{eqnarray}
Define
\begin{eqnarray}
A_i =\left(\frac{-\D K(a_i)}{K(a_i)}
\frac{1}{\frac{K(a_1)+K(a_2)}{2K(y_1)}-1}\right)^{1/2} \quad \mbox{ for  }
i=1,2.
\end{eqnarray}
Assume that
\begin{eqnarray*}
(H_1)\quad
\begin{cases}
 & \tilde c_2\geq \frac{K(a_1)+K(a_2)}{2K(y_1)}-1,\quad
-\D K(a_1)>0,\quad -\D K(a_2)>0\\
 & \forall\, \, (a_1,a_2)\in \nu^+(z_1)\times \nu^+(z_2) \mbox {
 such that } c_\infty (a_1,a_2)\leq c_\infty(y_1,y_1).
\end{cases}
\end{eqnarray*}
Then, on $\partial ([A_1,+\infty)\times [A_2,+\infty)\times \tilde
V)$, we have
\begin{eqnarray}
J(u)\geq c_\infty(y_1,y_1)
\end{eqnarray}
and $2)$ is therefore satisfied. We are now left with $3)$, that is
to verify \eqref{L'} for $(a_1,\l_1)$ and $(a_2,\l_2)$, $\l_1$ in
$(A_1,+\infty)$, $\l_2$ in $(A_2,+\infty)$. This amounts to
requiring, if we add the other requirement that $\l_i\rho_i$'s are
large,
\begin{eqnarray*}
(H_2)\quad
\begin{cases}
 & \frac{1}{\rho_i ^{n-4}A_i ^{n-6}}+\frac{|\n K(a_i)|^2}{
K(a_i)^2}+ \frac{|D^2 K(a_i)|^2}{A_i^2K(a_i)^2}+\sup_{B_i}\frac{|D^3 K|}{A_i K(a_i)}\leq c'_1\frac{-\D K(a_i)}{K(a_i)}\\
& (\sup K/ K(y_1))^{\max(1,(n-4)/8)}||f||K(y_1)^{(n-4)/8} \leq c_1''\\
& \sum (\sup K) K(y_1)^{-1}\log(\l_i\rho_i)^{(n+4)/n}(\l_i\rho_i)^{-4}\leq c_1''.\\
 & A_i\rho_i\geq \frac{1}{c'_1} ;\quad \rho_i \leq d(a_1,a_2)/3
\quad\forall \, \, i=1,2.
\end{cases}
\end{eqnarray*}
 Next we are going to show that $(H_2)$
follows from (for $C_0$, $C_1$ suitable small constants)
$$
(H_3)\,
\begin{cases}
 & w=w(a_1,a_2)=\frac{K(a_1)+K(a_2)}{2K(y_1)}-1\leq C_0,\\
 & w^{\frac{n-6}{n-4}}\bigl(\frac{1}{d(a_1,a_2)^2}+
\frac{1}{\rho_0^2}\bigr)+\frac{|\n K(a_i)|^2}{K(a_i)^2}
+w^{\frac{1}{3}}\sup_{B_i}\bigl(\frac{|D^3 K|}{ K(a_i)}\bigr)^{\frac{2}{3}}+ w^{\frac{1}{2}}\frac{|D^2K(a_i)|}{K(a_i)}\\
&\leq \frac{C_1}{ 1+ (\sup K/ K(y_1))^{\max(1,(n-4)/8)}}\frac{-\D K(a_i)}{K(a_i)}\\
 & \forall \, \, (a_1,a_2)\in \nu^+(z_1)\times \nu^+(z_2) \mbox{
 such that } c_\infty(a_1,a_2)\leq c_\infty(y_1,y_1),
\end{cases}
$$
where $\rho_0$ is any fixed positive constant 
Picking up any $\rho_0>0$, and choosing
\begin{eqnarray}
\tilde\rho_i =\min \biggl(\frac{d(a_1,a_2)}{3},
\rho_0\biggr),
\end{eqnarray}
 We now check that $A_i\tilde\rho_i\geq 1/c'_1$.
Indeed, using the first and the second conditions of \eqref{L'}, we obtain
\be
\left(A_i \tilde\rho_i\right)^2 \geq \frac{-\D K(a_i)}{9w K(a_i)} d(a_1,a_2)\geq C_1^{-1} w^{-2/(n-4)} \geq C_1^{-1}C_0^{-2/(n-4)}.
\ee
Since $C_1$ and $C_0$ are chosen small, this implies that $A_i\rho_i$ is very large. Notice that, by easy computations, the other conditions of $(H_2)$ follow from $(H_3)$

The fact that $\tau$ is 1 follows under \eqref{L'}. Using Theorem \ref{t:1}, we
derive the existence of a solution. The proof of Theorem \ref{t:2} is therefore
completed.
\end{pfn}

\vspace{3mm}\noindent
{\bf Acknowledgements.} Part of this work was done while I was visiting the Mathematics Department of the University of Roma ``La Sapienza''. I would like to thank the Mathematics Department for its warm hospitality. I also thank Professors A. Ambrosetti, A. Bahri, M. Grossi, F. Pacella and O. Rey for their encouragement and constant support over the years. I owe special thanks to Professor M. Ben Ayed for fruitful discussions and for his generous help during the preparation of this paper.

\end{document}